\documentclass[11pt]{article}
\usepackage{graphicx}
\usepackage{amsmath}
\usepackage{amssymb}
\usepackage{amsfonts}
\usepackage{color}
\usepackage{subfigure}
 \usepackage{cases}
\usepackage{setspace}
\usepackage{fullpage}
\usepackage{listings}

\newtheorem{definition}{Definition}
\newtheorem{remark}{Remark}

\newcommand{\nc}{\newcommand}
\nc{\R}{\Bbb{R}}
\nc{\Z}{\Bbb{Z}}
\nc{\Pp}{\Bbb{P}}
\nc{\Ap}{\Bbb{A}}
\nc{\Wp}{\Bbb{W}}
 \nc{\brho}{\boldsymbol \rho}
\nc{\va}{\vec{\boldsymbol \mu}}
\nc{\ve}{\vec{\boldsymbol \epsilon}}
\nc{\bk}{{\bf k}}
\nc{\vrho}{\vec{\brho}} 
\nc{\vr}{\vec{\bf r}}
\nc{\bx}{{\bf x}}
\nc{\vx}{\vec{\bf x}}
\nc{\om}{\omega}
\nc{\brhoi}{\brho^{\cal I}}
\nc{\bzetai}{\bzeta^{\cal I}}
\nc{\vzeta}{\vec{\bzeta}}
\nc{\vzetai}{\vec{\bzeta}^{\cal I}}
\nc{\cI}{{\cal I}}
\nc{\cJ}{{\cal J}}
\nc{\cF}{{\cal F}}
\nc{\cW}{{\cal W}}
\nc{\cA}{{\cal A}}
\nc{\cL}{{\cal L}}
\nc{\cS}{{\cal S}}
\nc{\cC}{{\cal C}}
\nc{\cN}{{\cal N}}
\nc{\cM}{{\cal M}}
\nc{\vrhoi}{\vec{\brho}^{\,\cal I}}
\nc{\xii}{\xi^{\cI}}
\nc{\etai}{\eta^{\cI}}
\nc{\la}{\lambda}
\nc{\de}{\delta}
\nc{\ep}{\varepsilon}
\nc{\vu}{\vec{\bf u}}
\nc{\bu}{{\bf u}}
\nc{\vui}{\vec{\bf u}^{\cI}}
\nc{\bui}{{\bf u}^{\cI}}
\nc{\bt}{{\bf t}}
\nc{\vt}{\vec{\bt}}
\nc{\bn}{{\bf n}}
\nc{\vn}{\vec{\bn}}
\nc{\bm}{{\bf m}}
\nc{\vm}{\vec{\bm}}
\nc{\vrp}{\vec{{\bf r}'}}
\nc{\vrc}{\vec{{\bf r}^c_p}}
\nc{\ts}{\tilde s}
\nc{\os}{\overline s}
\nc{\tom}{\tilde \om}
\nc{\tO}{\tilde \Omega}
\nc{\tS}{\tilde S}
\nc{\oS}{\overline S}
\nc{\vrhos}{\vrho_{\star}}
\nc{\vrhosi}{\vrho_{\star}^{\cI}}
\nc{\brhos}{\brho_{\star}}
\nc{\brhosi}{\brhos^{\cI}}
\nc{\vms}{\vm_{\star}}
\nc{\vmi}{\vm^{_{\cI}}}
\nc{\vM}{\vec{\bf M}}
\nc{\vMi}{\vM_{_{\cI}}}
\nc{\Ppi}{\Pp_{\cI}}
\nc{\vxi}{\vec{\bxi}}
\nc{\bK}{{\bf K}}
\nc{\bmi}{\bm_{_{\cI}}}
\nc{\Pppi}{\Ppi^p}
\nc{\be}{{\bf e}}
\nc{\bep}{{\bf e}^p}
\renewcommand{\hat}{\widehat}

\begin{document}
\title{Motion Estimation and Imaging of Complex Scenes\\ with
Synthetic Aperture Radar } \author{Liliana Borcea\footnotemark[2],
Thomas Callaghan\footnotemark[2], and George
Papanicolaou\footnotemark[3]
}\renewcommand{\thefootnote}{\fnsymbol{footnote}}
\footnotetext[2]{Computational and Applied Mathematics, Rice
University, MS 134, Houston, TX 77005-1892. (borcea@caam.rice.edu and
tscallaghan@rice.edu)} \footnotetext[3]{Department of Mathematics,
Stanford University, Stanford, CA 94305. (papanicolaou@stanford.edu)}
\date{} \maketitle
%  \tableofcontents
\numberwithin{equation}{section} %sets equation numbers <chapter>.<section>.<index>

% ------------------------------
\begin{abstract}
We study synthetic aperture radar (SAR) imaging and motion estimation
of complex scenes consisting of stationary and moving targets.  We use
the classic SAR setup with a single antenna emitting signals and
receiving the echoes from the scene. The known motion estimation
methods for such setups work only in simple cases, with one or a few
targets in the same motion. We propose to extend the applicability of
these methods to complex scenes, by complementing them with a data
pre-processing step intended to separate the echoes from the
stationary targets and the moving ones. We present two approaches.
The first is an iteration designed to subtract the echoes from the
stationary targets one by one. It estimates the location of each
stationary target from a preliminary image, and then uses it to define
a filter that removes its echo from the data. The second approach is
based on the robust principle component analysis (PCA) method. The key
observation is that with appropriate pre-processing and windowing, the
discrete samples of the stationary target echoes form a low rank
matrix, whereas the samples of a few moving target echoes form a high
rank sparse matrix. The robust PCA method is designed to separate the
low rank from the sparse part, and thus can be used for the SAR data
separation.  We present a brief analysis of the two methods and
explain how they can be combined to improve the data separation for
extended and complex imaging scenes.  We also assess the performance
of the methods with extensive numerical simulations.

\end{abstract}
% ------------------------------
\section{Introduction}
\label{sect:intro}
\setcounter{equation}{0}

\begin{figure}[t]
\vspace{-0.4in}
  \begin{center}
    \input{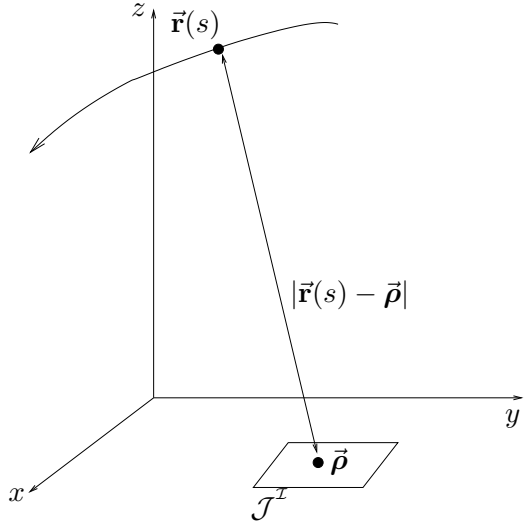}
  \end{center}
  \caption{Setup for synthetic aperture imaging.}
  \label{fig:setup}
\end{figure}

Synthetic aperture radar (SAR) is an important technology capable of
computing high resolution images of a scene $\cJ^{^\cI}$ on the ground,
using measurements made from an antenna mounted on a platform
flying over it \cite{cheney2001mathematical,Jakowatz, Curlander}.  The
antenna periodically emits a signal $f(t)$ and records the echoes, the
data that are processed to form an image.  The problem setup is
illustrated in Figure \ref{fig:setup}.  The data $D(s,t)$ are indexed
by the {\em slow-time} $s$ and the {\em fast-time} $t$.  The slow-time
parameterizes the antenna position $\vr(s)$ at the moment it emits the
signal. The fast-time parameterizes the time between consecutive
illuminations at $s$ and $s+\Delta s$, so that $0<t<\Delta s.$

The recordings $D(s,t)$ are approximately and up to a multiplicative
factor a linear superposition of the probing signals $f(t)$
time-delayed by the round-trip travel time between the antenna
position $\vr(s)$ and the locations $\vrho$ of the scatterers
(targets) in  $\cJ^{^\cI}$. Assuming that the waves propagate at the
constant speed of light $c$, the round-trip travel time is given by
\begin{equation}
\tau(s,\vrho)=2|\vr(s)-\vrho|/c.
\end{equation}
To form an image, the data are first convolved with the time reversed
emitted signal $f(t)$, a process known as match-filtering or pulse
compression. Then, they are back-propagated to search points $\vrhoi
\in \cJ^{^\cI}$ using the travel times $\tau(s,\vrhoi)$. The image is
obtained by superposing the pulse compressed and back-propagated data
over a slow time interval defining the synthetic aperture $a$, the
distance traveled by the antenna.  

SAR systems that cover long apertures and emit broad-band high
frequency signals give very good resolution of images of stationary
scenes.  The resolution can be of the order of ten centimeters at
ranges of ten kilometers away. However, the images are degraded when
there is motion in the scene.  Depending on how fast the targets move,
they may appear blurred and displaced, or they may not be visible at
all. The motion needs to be estimated, and then compensated in the
image formation, in order to bring these targets in focus.  This is a
complicated task for complex imaging scenes, consisting of many
stationary targets and moving ones.

There are basically two approaches to motion estimation. The first
determines the motion from the phase modulations of the return echoes
\cite{sar,ding2000analysis,ding2002time,barbarossa1992detection,
  sparr-time, zhuMTI, fienup, jao, kirscht, perry,ender1993}, assuming
either a single target in the scene, or many targets moving the same
way. It does not work well in complex scenes. The estimation is
sensitive to the presence of strong stationary targets and to targets
that have different motion. The second approach uses multiple
receiver and/or transmitter antennas
\cite{friedlander,wang2004,wang2006}. It forms a collection of images
with the echoes measured by each receiver-transmitter pair. Then,
it extracts the target velocities from the phase variation of the images
with respect to the receiver/transmitter offsets.

We are interested in the first approach that uses the classic
SAR setup with a single antenna. To extend its applicability to
complex scenes, we complement it with a data pre-processing step
intended to separate the echoes from the stationary targets and the
moving ones.  A successful separation allows motion estimation to
be carried with the echoes from the moving targets alone, using the
algorithms in
\cite{sar,ding2000analysis,ding2002time,barbarossa1992detection,
  sparr-time, zhuMTI, fienup, jao, kirscht, perry,ender1993}.

We present two data separation methods.  The first extends the ideas
in \cite{delcueto,delcueto2}, developed in a different context of
imaging in scattering layered media.  The method seeks to subtract the
stationary target echoes one by one, using a combination of travel
time transformations and a linear annihilation filter.  Each travel
time transformation is relative to one stationary target at a time,
whose location can be estimated from a preliminary image. After the
transformation, the echo from this target becomes essentially
independent of the slow time, and so it can be annihilated with a linear
filter such as a difference operator in $s$.  The travel time
transformation is undone after the filtering, and the method moves on
to the next stationary target.

The second method is based on the robust principle component analysis
(PCA) method \cite{candesRPCA}. We have shown in \cite{BCP12} with
analysis and numerical simulations that after appropriate
pre-processing and windowing, the samples of the echoes from
stationary targets form a low rank matrix. Moreover, the samples of
the echoes from a few moving targets form a high rank but sparse
matrix. Thus, the echoes can be separated with the robust PCA method,
which is designed to decompose a matrix into its low rank and sparse
parts \cite{candesRPCA}.  We review briefly the results in
\cite{BCP12}. We also explain that in practice it may be useful to use
a combination of the two approaches in order to achieve better data
separation for extended and complex imaging scenes.

The paper is organized as follows.  We begin in section
\ref{sect:prelims} with a brief description of basic SAR data
processing and image formation.  The data separation methods are
described in section \ref{sect:algdesc}. We analyze them in section
\ref{sect:analysis} for simple imaging scenes. Then we show with
numerical simulations in section \ref{sect:numeric} that the
methods can handle complex scenes. We end with a summary in section
\ref{sect:conc}.
% ------------------------------------------
\section{Preliminaries}
\setcounter{equation}{0}
\label{sect:prelims}

We begin with two data pre-processing steps: pulse compression and
range compression.  They are commonly used in SAR, but in addition,
they are an important first step in our data separation methods. We
also describe the generic SAR image formation. The goal of the paper
is to explain how to complement the image formation with data
separation and filtering in order to estimate motion and focus images
of complex scenes.

% ------------------------------------------
\subsection{Pulse and range compression}
\label{sect:compress}

Typically, the antenna emits signals that consist of a base-band
waveform $f_B(t)$ modulated by a carrier frequency
$\nu_0=\omega_o/(2\pi)$
\begin{equation}
f(t)=\cos(\omega_ot)f_B(t).
\label{eq:signal}
\end{equation}
The Fourier transform of the signal is
\begin{equation}
\hat{f}(\omega)=\int_{-\infty}^\infty dt \, f(t)e^{i\omega
  t}=\frac{1}{2}\left[\hat f_B(\omega+\omega_o)+\hat
  f_B(\omega-\omega_o)\right],
\end{equation}
where the support of $\hat f_B(\omega)$ is the interval $[-\pi B,\pi
B]$, and $B$ is the bandwidth.  Ideally, the waveform $f_B(t)$ should
be a pulse of short duration, so that travel times can be easily
estimated. It should also carry enough power so that the received
echoes are stronger than the ambient noise.  However, antennas have
instantaneous power limitations and they can transmit large net power only 
by spreading it over long signals, like chirps.  Consequently, the
measured echoes are long and it is impossible to read travel times
directly from them. This difficulty is overcome by compressing the
echoes, as if they were due to short emitted pulses.  

Pulse compression is achieved by convolving (match filtering) the data
with the complex conjugate of the time reversed emitted signal
\begin{equation}
D_p(s,t)=\int dt' D(s,t')\overline{f(t'-t)}.
\end{equation}
Mathematically, this is equivalent to the antenna emitting the signal
\begin{equation}
f_p(t)=\int dt' f(t')\overline{f(t'-t)},
\end{equation}
and receiving the echoes $D_p(s,t)$. The signal $f_p(t)$ is a pulse of
duration of the order $1/B$, as shown for example in \cite{Jakowatz}.

Range compression is another important data pre-processing step.  It
consists of shifting the fast time $t$ in the argument of $D_p(s,t)$
by the travel time to a reference point $\vrho_o\in\cJ^{^\cI}$,
\begin{equation}
D_r(s,t')=D_p(s,t'+\tau(s,\vrho_o)).
\label{eq:rangecompress}
\end{equation}
Here $t'$ is the shifted fast-time, satisfying
$t'+\tau(s,\vrho_o)\in[0,\Delta s].  $ 

\begin{figure}[t]
\centering
\includegraphics[width=.5\columnwidth]{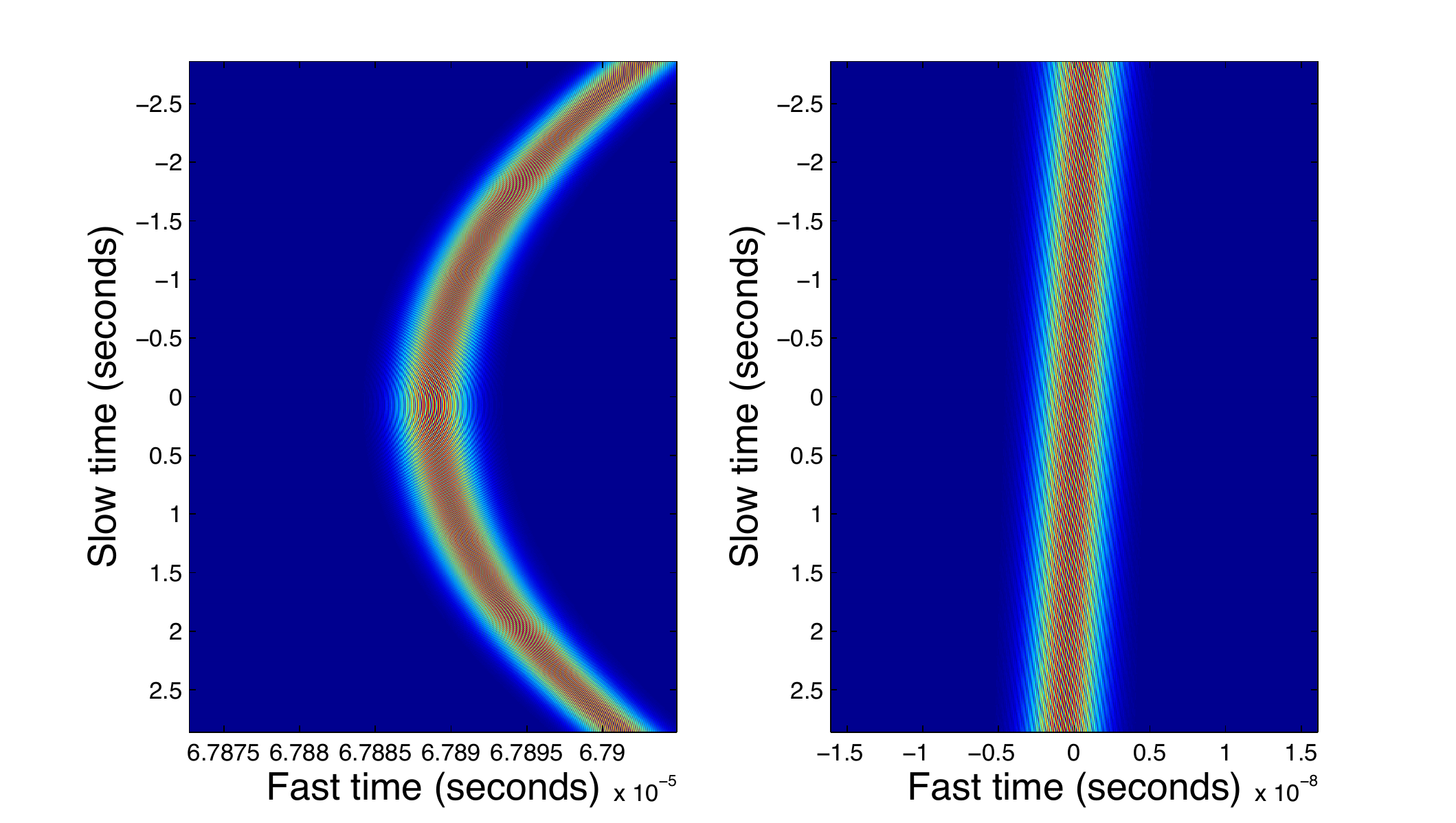}
\caption{Left: Pulse compressed data $D_p(s,t)=f(t-\tau(s,\vrho))$
  from a target at $\vrho=(0,5,0)$m, recorded by a SAR antenna moving
  at speed $V=70$ m/s along a linear aperture of $400$ m, at range $L
  = 10$ km from $\vrho$.  Right: Range-compressed data $D_r(s,t')$
  with respect to $\vrho_o=(0,0,0)m$.  }
\label{fig:hyperbola}
\end{figure}

Range compression is beneficial for numerical computations and it is
essential for our data separation methods because it removes partially
the dependence of $D_p(s,t)$ with respect to the slow time. We
illustrate this in Figure \ref{fig:hyperbola}, where we compare the
pulse compressed echo from a stationary target before and after the
range compression.
The target is at location $\vrho=(0,5,0)$m and the reference point is
at $\vrho_o = (0,0,0)$m.  The SAR platform moves at speed $70$m/s
along a linear aperture $a = 400$m at range $L=10$ km from
$\vrho$. The imaging scene lies on a flat surface at zero altitude.
We show in the left plot in Figure \ref{fig:hyperbola} the amplitude
of the pulse compressed echo as a function of $s$ and $t$. Note that
its peak lies on the curve defined by
\[
\left(\frac{ct}{2}\right)^2=|\vr(s)-\vrho|^2.
\]
This curve is a hyperbola because the aperture is linear and the
platform moves at constant speed. The right plot in Figure
\ref{fig:hyperbola} shows the amplitude of the range and pulse
compressed echo. Its peak lies on the curve 
\[
\frac{ct'}{2}=|\vr(s)-\vrho|-|\vr(s)-\vrho_o|,
\]
which looks almost as a vertical line segment, because the dependence
on the slow time has been approximately removed by the range
compression. 

We work henceforth with the range and pulse compressed data and drop
the prime from the shifted fast time $t'$ to simplify the notation.
Borrowing terminology from geophysics, we call $D_r(s,t)$ the {\em
data traces} or simply the {\em traces}.

% ------------------------------------------
\subsection{Image formation}
\label{sect:image}
A SAR image is formed by superposing the traces over the aperture and
back propagating them to the imaging points $\vrhoi$ using travel
times
\begin{align}
\label{eq:image}
\cI(\vrhoi)&=\sum_{j=-n/2}^{n/2}
D_r\left(s_j,\tau\left(s_j,\vrhoi\right)-\tau(s_j,\vrho_o)\right). 
% \\ &\approx \frac{1}{\Delta s}\int_{-S(a)}^{S(a)} ds.
% D_r\left(s,\tau\left(s,\vrhoi\right)-\tau(s,\vrho_o)\right)\nonumber.
\end{align}
The slow time samples $s_j$ discretize the interval $[-S(a),S(a)]$
which defines the aperture $a$ along the trajectory of the SAR
platform.  The sampling $\Delta s$ is uniform and equals the time
interval between tho consecutive signal emissions from the antenna.
To simplify the notation, we assume that the origin of the slow time
corresponds to the center of the aperture. Thus, $n$ is an even
integer and
\[
2 S(a) = n \Delta s.
\]

Equation (\ref{eq:image}) is a simplification of the imaging function
used in practice, which includes a weight factor that accounts for
geometrical spreading. The weight makes a difference only when the
aperture is large. We neglect it in this paper because we focus
attention on motion estimation, which is necessarily done over small
successive sub-apertures. Targets may have complicated trajectories
over long data acquisition times, but we can approximate their motion
by a uniform translation over short durations $2 S(a)$ defining
small sub-apertures $a$. The idea is that the speed of translation can
be estimated sub-aperture by sub-aperture, and then it can be
compensated in the image formation. By compensation we mean that if we
seek to image a moving target and we estimate its speed to be
$\vec{\bf u}$, we can form an image as
\begin{equation}
\label{eq:imagecompu}
\cI_{_{\vec{\bf u}}}(\vrhoi)=\sum_{j=-n/2}^{n/2}
D_r\left(s_j,\tau\left(s_j,\vrhoi + s_j \vec{\bf
u}\right)-\tau(s_j,\vrho_o)\right).
\end{equation}
We repeat the process for successive sub-apertures, and assemble the
final image by superposing the results. This superposition requires
proper weighting to compensate geometrical spreading effects over long
SAR platform trajectories.

Motion estimation and its compensation in imaging works for simple
scenes with one or a few targets that move at approximately the same
speed. To our knowledge, there is no approach that estimates motion in
complex scenes with a single antenna. Moreover, even if a target
speed were estimated, the compensation in (\ref{eq:imagecompu}) would
bring that target in focus but it would blur the others.  To address
the difficulty of imaging with motion estimation in complex scenes, we
propose a data pre-processing step for separating the traces from the
stationary targets and the moving ones. If we could achieve such
separation, we could carry the motion estimation on the data traces
from the moving targets alone. Moreover, we could image the stationary
and moving targets separately, and then combine the results to obtain
an image of the complex scene.
% ------------------------------------------
\section{Data separation}
\label{sect:algdesc}
\setcounter{equation}{0}

We present two complementary approaches for separating traces from
stationary and moving targets.  They are analyzed in section
\ref{sect:analysis}, in simple setups.  Here we describe how they work
and illustrate the ideas with numerical simulations, for a complex
scene named hereafter ``Scene 1''.  It consists of twenty stationary
targets and two moving targets with speeds of $28$m/s and $14 $m/s,
respectively. We refer to section \ref{sect:numeric} for details on
the setup of the numerical simulations.  The traces are displayed in
the left plot of Figure \ref{fig:dataDecomp}.
% \begin{figure}[!t]
% \centering\includegraphics[width=.6\columnwidth]{./Figures/data2mov}
% \caption{Traces for a scene consisting of twenty stationary targets
%  and two moving targets with speed $28$m/s and $14$m/s respectively.
%  The vertical traces correspond to the stationary targets, while the
%  sloped traces correspond to the moving targets.}
% \label{fig:data}
% \end{figure}
% ------------------------------------------
\subsection{Travel time transformations and data separation}
\label{sect:ANNIH}
\begin{figure}[!t]
\centering\includegraphics[width=.85\columnwidth]{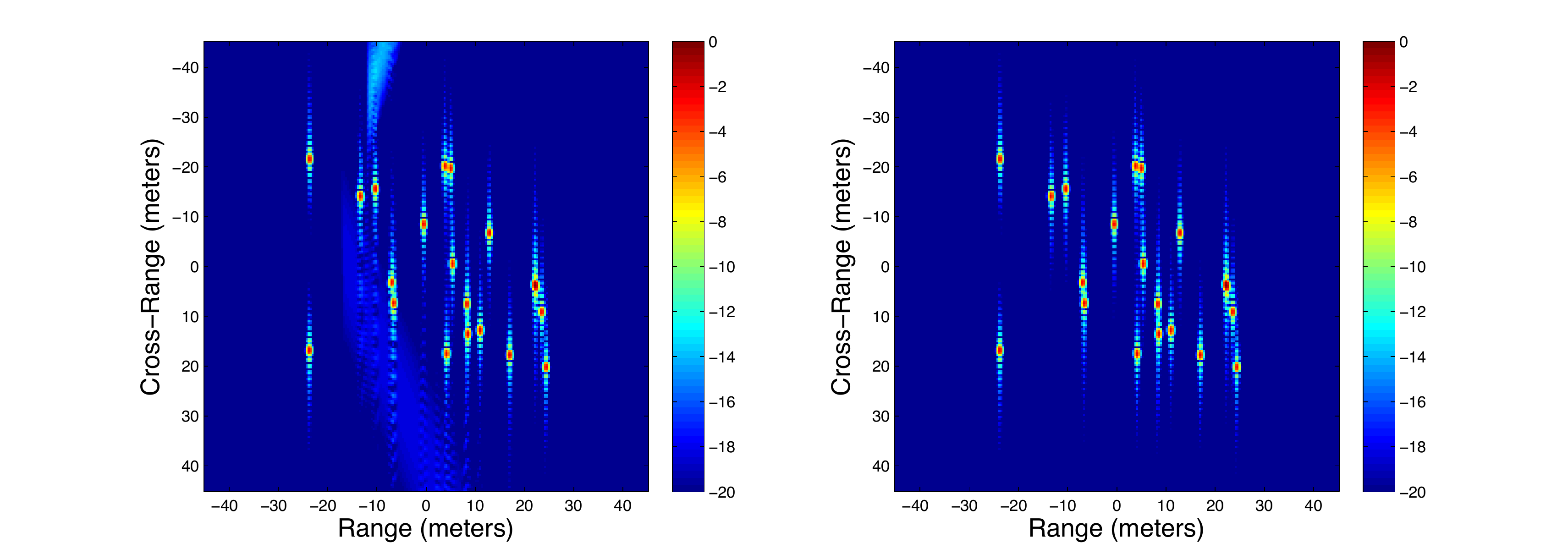}
\caption{Left: The image of the complex Scene 1, with twenty
  stationary targets and two moving targets.  Right: Image of the same
  stationary scene, and no moving targets. The effect of the moving
  targets in the left image amounts to two faint streaks.}
\label{fig:imageComp}
\end{figure}

The basic idea of the first data separation approach is that if we
have an estimate $\vrho^e$ of the location $\vrho$ of a stationary
target, we can shift the time $t$ by the travel time $\tau(s,\vrho^e)$
to remove approximately the dependence on the slow time $s$ of the
trace from that target. Then, we can subtract the trace from the data
by exploiting its independence of $s$.  To determine $\vrho^e$, we
need a preliminary image of the scene. The assumption is that the
stationary targets dominate the scene, and therefore they are visible
in the preliminary image even in the presence of moving targets. See
Figure \ref{fig:imageComp} for an illustration.

\begin{figure}[!t]
\centering\includegraphics[width=.5\columnwidth]{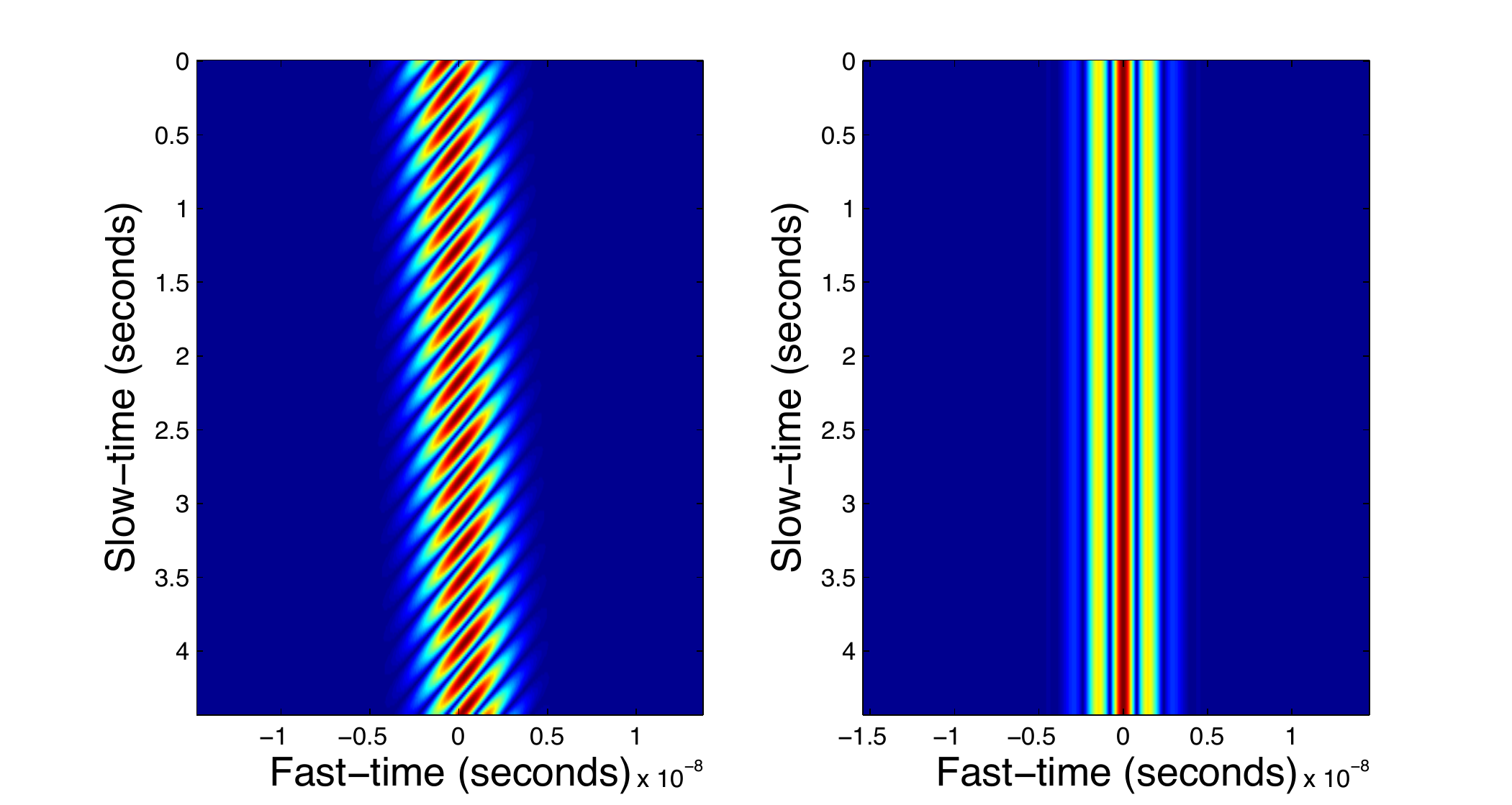}
\caption{Left: The trace $|D_r(s,t)|$ from one stationary target.
  Right: the traces after the travel time transformation relative to
  the exact location of the target. Note that there is no dependence
  on $s$ after the transformation.}
\label{fig:straightTrace}
\end{figure}

Let $\mathbb{T}_+^{\vrho^e}$ be the map taking the range compressed
data $D_r(s,t)$ to the data in the transformed coordinates
$(s,t+\Delta \tau(s,\vrho^e))$
\begin{equation}
\label{eq:TT+}
[\mathbb T_+^{\vrho^e} D_r](s,t)=D_r(s,t+\Delta\tau(s,\vrho^e)),
\end{equation}
where 
\begin{equation}
\label{eq:DeltaTau}
\Delta \tau(s,\vrho^e) = \tau(s,\vrho^e) - \tau(s,\vrho_o).
\end{equation}
Recall that $\vrho_o$ is the reference point in the range compression.
To illustrate the effect of the map (\ref{eq:TT+}), we show on the
left in Figure \ref{fig:straightTrace} the trace from a target at
location $\vrho \ne \vrho_o$. The plot on the right shows the trace
after the travel time transformation (\ref{eq:TT+}), with the ideal
estimate $\vrho^e = \vrho$. The transformation makes the trace
independent of $s$, and thus it can be annihilated by a difference
operator in $s$, 
\begin{equation}
\label{eq:DiffOp}
\mathbb{D}_s[\mathbb T_{+}^{\vrho^e} D_r](s,t) = \frac{1}{\Delta s}
\left\{[\mathbb T_{+}^{\vrho^e} D_r](s+\Delta s,t) - [\mathbb
  T_{+}^{\vrho^e} D_r](s,t)\right\} \approx \frac{\partial}{\partial
  s} [\mathbb T_{+}^{\vrho^e} D_r](s,t).
\end{equation}
This is an extension of the ideas proposed in
\cite{delcueto,delcueto2} for annihilating echoes from a layered
scattering medium. Other annihilation operators than (\ref{eq:DiffOp})
may be used, as well. In practice, where noise plays a role, they
should be complemented with a smoothing process. We refer to
\cite{delcueto,delcueto2} for a detailed analysis of annihilator
filters. For the simulations in this paper the annihilator
$\mathbb{D}_s$ works well, and there is no noise added to the data.

\begin{figure}[!t]
\centering\includegraphics[width=.7\columnwidth]{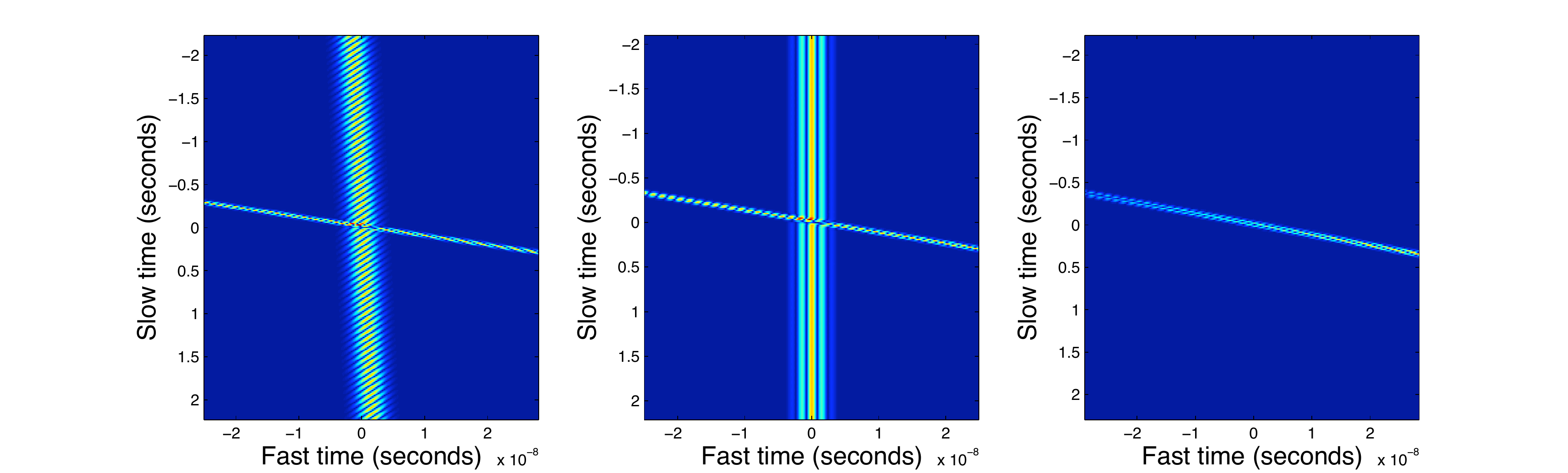}
\caption{Left: The traces $D_r(s,t)$ for a scene with one stationary
  target at location $\vrho$ and a moving one at $28$m/s.  Middle: The
  traces after the travel time transformation $[\mathbb T_+^{\vrho}
    D_r](s,t)$.  Right: The filtered traces $[ \mathcal
    Q^{\vrho^e}D_r](s,t)$.}
\label{fig:rotScatt}
\end{figure}

After taking the derivative, we map the traces back to the $(s,t)$
coordinates using the travel time transformation
$\mathbb{T}_{-}^{\vrho^e}$. The data filter is a linear operator
denoted by $Q^{\vrho^e}$, given by composition of
$\mathbb{T}_{+}^{\vrho^e}$, $\mathbb{D}_s$ and
$\mathbb{T}_{-}^{\vrho^e}$,
\begin{eqnarray}
\label{eq:OP1}
\left[\mathcal Q^{\vrho^e} D_r\right](s,t)= \left[\mathbb
T_{-}^{\vrho^e}\mathbb{D}_s \mathbb T_+^{\vrho^e} D_r\right](s,t)
= \left[\mathbb{D}_s
D_r(s,t'+\Delta\tau(s,\vrho^e))\right]_{t'=t-\Delta\tau(s,\vrho^e)}.
\end{eqnarray}
We demonstrate in Figure \ref{fig:rotScatt} the travel-time
transformations and annihilation steps for a scene with one stationary
and one moving target. Note that the trace from the stationary target
located at $\vrho$ is completely removed by the mapping $Q^{\vrho^e}$,
because $\vrho^e = \vrho$.

If there are $N$ stationary targets in the scene, we can subtract
their traces from the data one by one, by iterating the procedure
above. Let $\vrho_j^e$ be the estimated locations of the targets, with
$j = 1, \ldots, N$. The data filter $\mathcal Q$ is the composition of
the linear operators (\ref{eq:OP1}),
\begin{equation}
\label{eq:OP2}
\left[\mathcal QD_r\right](s,t)=\left[\mathcal Q^{\vrho_N^e}
  \circ\mathcal Q^{\vrho_{N-1}^e} \circ \ldots \circ \mathcal
  Q^{\vrho_1^e} D_r\right] (s,t).
\end{equation}

The estimates $\vrho_j^e$ are not exact in practice, but the filter
(\ref{eq:OP1}) is robust with respect to the estimation errors. We
show this in section \ref{sect:analysis}, where we present a brief
analysis of the filters. In fact, our simulations show that it is
typical that one application of $\mathcal Q^{\vrho^e}$ removes at once
the traces from all the stationary targets located near
$\vrho^e$. Thus, the needed number of iterations in (\ref{eq:OP2}) may
be much smaller than the number of stationary targets in practice.

Because antennas measure discrete samples of the data $D_r(s,t)$, we
need to use interpolation when making the travel-time transformations.
We do so by implementing the travel-time transformations as phase
modulations in the frequency domain
\begin{equation}
\label{eq:FFT}
\left[\mathbb T_+^{\vrho^e} D_r\right](s,t)=\frac{1}{2\pi}\int
d\omega\ \hat D_r(s,\omega)e^{-i\omega(t+\Delta\tau(s,\vrho))},
\end{equation}
where
\[
\hat D_r(s,\omega)=\int dt\ e^{i\omega t} D_r(s,t).
\]
We implement (\ref{eq:FFT}) in the discrete setting using the fast
Fourier transform.

% ------------------------------------------
\subsection{Data separation with robust PCA}
\label{sect:RPCA}
We introduced recently in \cite{BCP12} a method for separating the
traces from stationary and moving targets using robust principal
component analysis (PCA). We review here the data separation method,
and give some details of its analysis in section \ref{sect:analysis}.

The robust PCA method was introduced and analyzed in
\cite{candesRPCA}.  It decomposes a matrix $\cM$ into a low rank matrix
$\cL$ and a sparse matrix $\cS$.  In our context, the entries in $\cM$
are the samples of the traces
\begin{equation}
\cM_{j\ell}=D_r(s_j,t_\ell),
\end{equation}
where 
\begin{equation}
s_j=j\Delta s,\quad j=-n/2,\ldots,n/2, \quad 2S(a) = n \Delta s,
\label{eq:discreteST}
\end{equation}
and
\begin{equation}
t_\ell=\ell\Delta s,\quad \ell=-m/2,\ldots,m/2, \quad \Delta s=m\Delta
t.
\label{eq:discreteFT}
\end{equation}
The sampling is at uniform intervals $\Delta s$ and $\Delta t$,
determined in practice by the SAR antenna system. To simplify the
notation, we take the origin of $s$ at the center of the aperture,
and the origin of $t$ in the middle of the interval $\Delta s$ between
two consecutive signal emissions. 

\begin{figure}[!t]
\centering\includegraphics[width=.75\columnwidth]{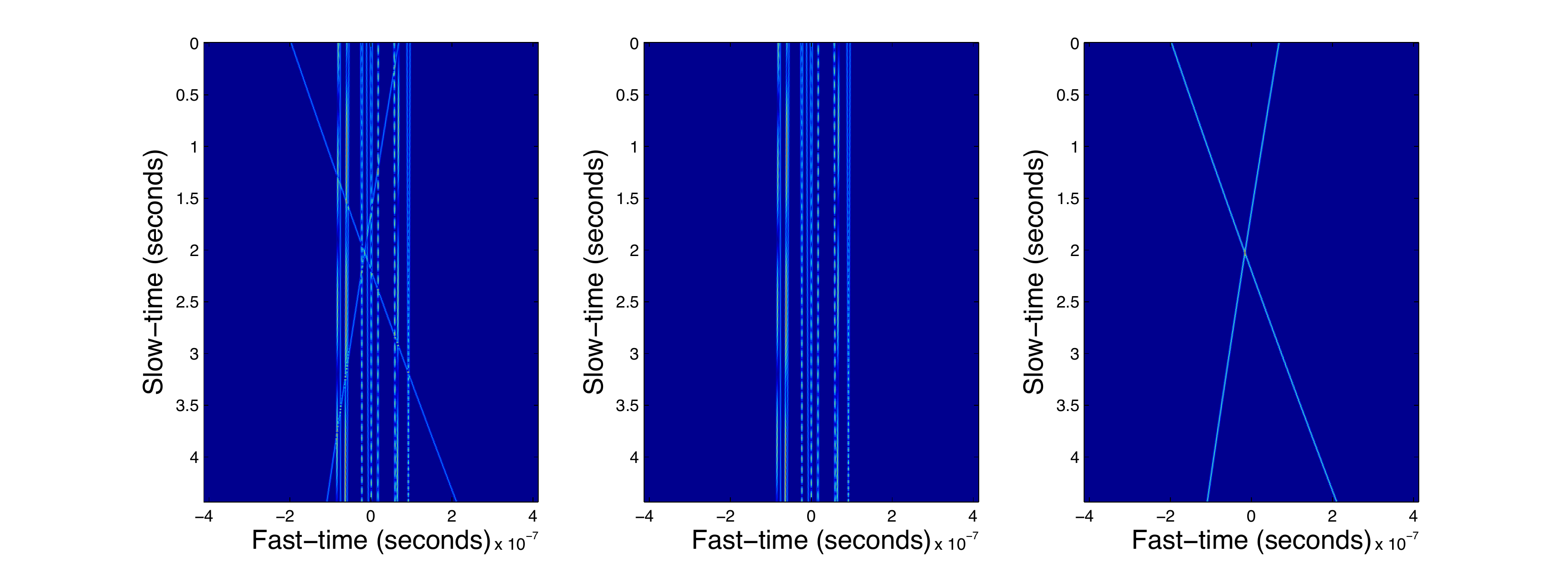}
\caption{Left: The matrix $\cM$ of the sampled traces from Scene 1.
  Middle: The matrix $\cL$ consisting of the traces from the
  stationary targets in Scene 1.  Right: The matrix $\cS$ consisting
  of the traces from the moving targets in Scene 1.}
\label{fig:dataDecomp}
\end{figure}
To explain why we can use robust PCA for the data separation, consider
the matrix $\cM$ of sampled traces for Scene 1. We display it in the
left plot of Figure \ref{fig:dataDecomp}. The traces that appear
almost vertical correspond to the stationary targets. This is because
these targets are not too far from the reference point $\vrho_o$, and
the range compression removes most of their $s$ dependence. We expect
that these traces form a low rank matrix $\cL$. We display them in the
middle plot of Figure \ref{fig:dataDecomp}. The range compression does
not remove the slow time dependence of the traces from the moving
targets, so they appear sloped in Figure \ref{fig:dataDecomp}.  They
form a high rank but sparse matrix $\cS$, shown in the right plot of
Figure \ref{fig:dataDecomp}. 

The robust PCA method is designed to separate $\cL$ from $\cS$ in $ \cM
= \cL + \cS.  $ It does so by solving an optimization problem called
principle component pursuit \cite{candesRPCA}
\begin{align} 
\label{eq:OPTIMIZATION}
\min_{\cL,\cS\in\mathbb
R^{n+1\times m+1}} \|\cL\|_*+\eta\|\cS\|_1, \quad \mbox{subject to}\quad
\cL+\cS=\cM. \end{align}
Here $\|\cL\|_*$ is the nuclear norm of $\cL$, the
sum of the singular values of $\cL$, $\|\cS\|_1$ is the matrix
1-norm of $\cS$ and 
\begin{equation}
\eta=\frac{1}{\sqrt{\max\{n+1,m+1\}}}.  
\end{equation} 
The optimization (\ref{eq:OPTIMIZATION}) can be applied to any matrix
$\cM$. But it is shown in \cite{candesRPCA} that if $\cM$ is a matrix with
a low rank part $\cL$ and a high rank sparse part $\cS$, then the
algorithm recovers exactly $\cL$ and $\cS$. We refer to
\cite{candesRPCA} for sufficient (not necessary) conditions
satisfied by $\cM$, so that it can be succesfully decomposed by
principle component pursuit.

Robust PCA is well suited for our purpose, but it cannot be applied as
a black box to the matrix of traces. We need to complement it with
properly calibrated data windowing.  To illustrate why, we present two
examples.  

\begin{figure}[t]
\begin{center}
 \subfigure{\includegraphics[width=.3\columnwidth]{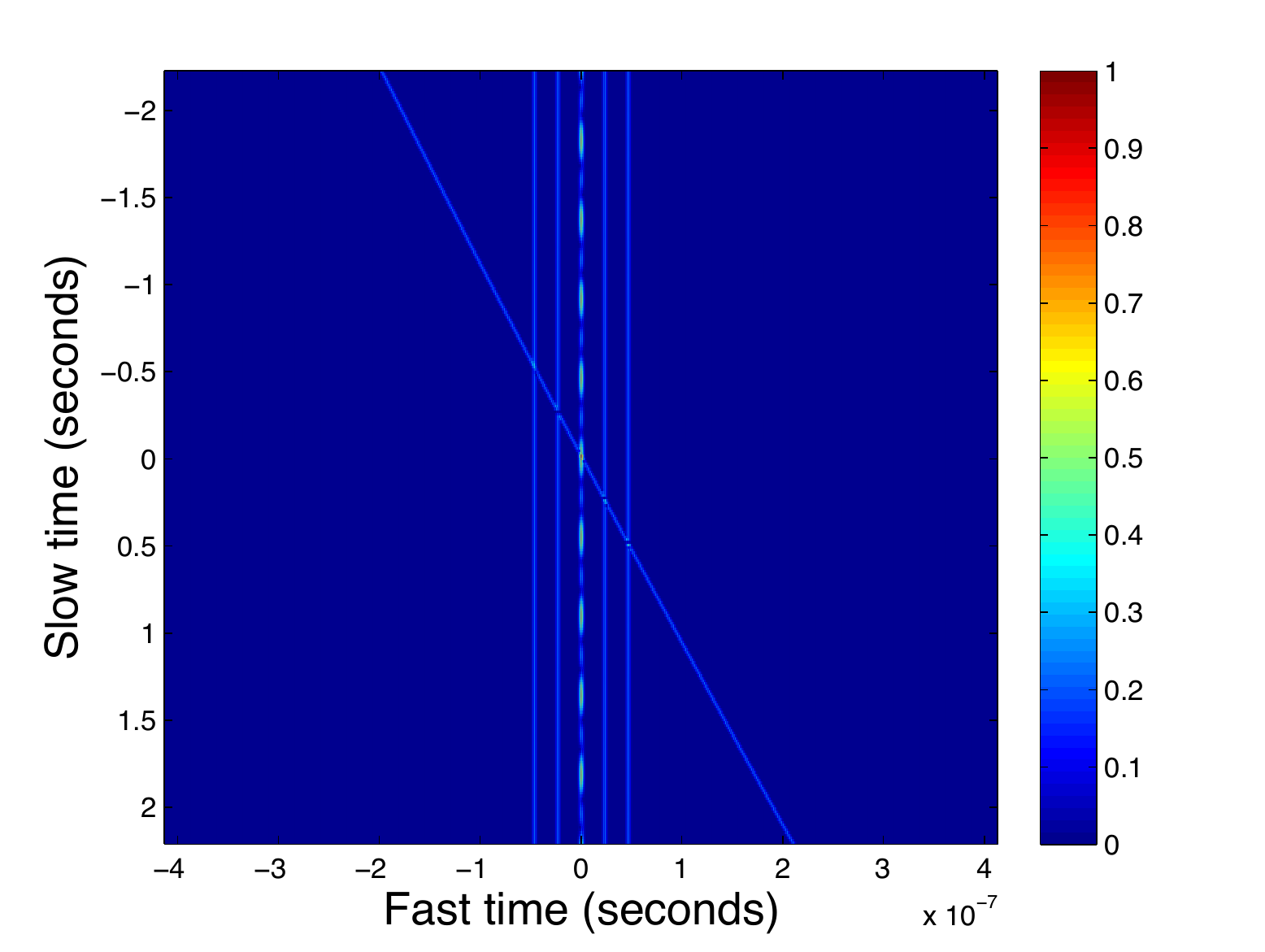}}\hspace{-0.2in}
 \subfigure{\includegraphics[width=.3\columnwidth]{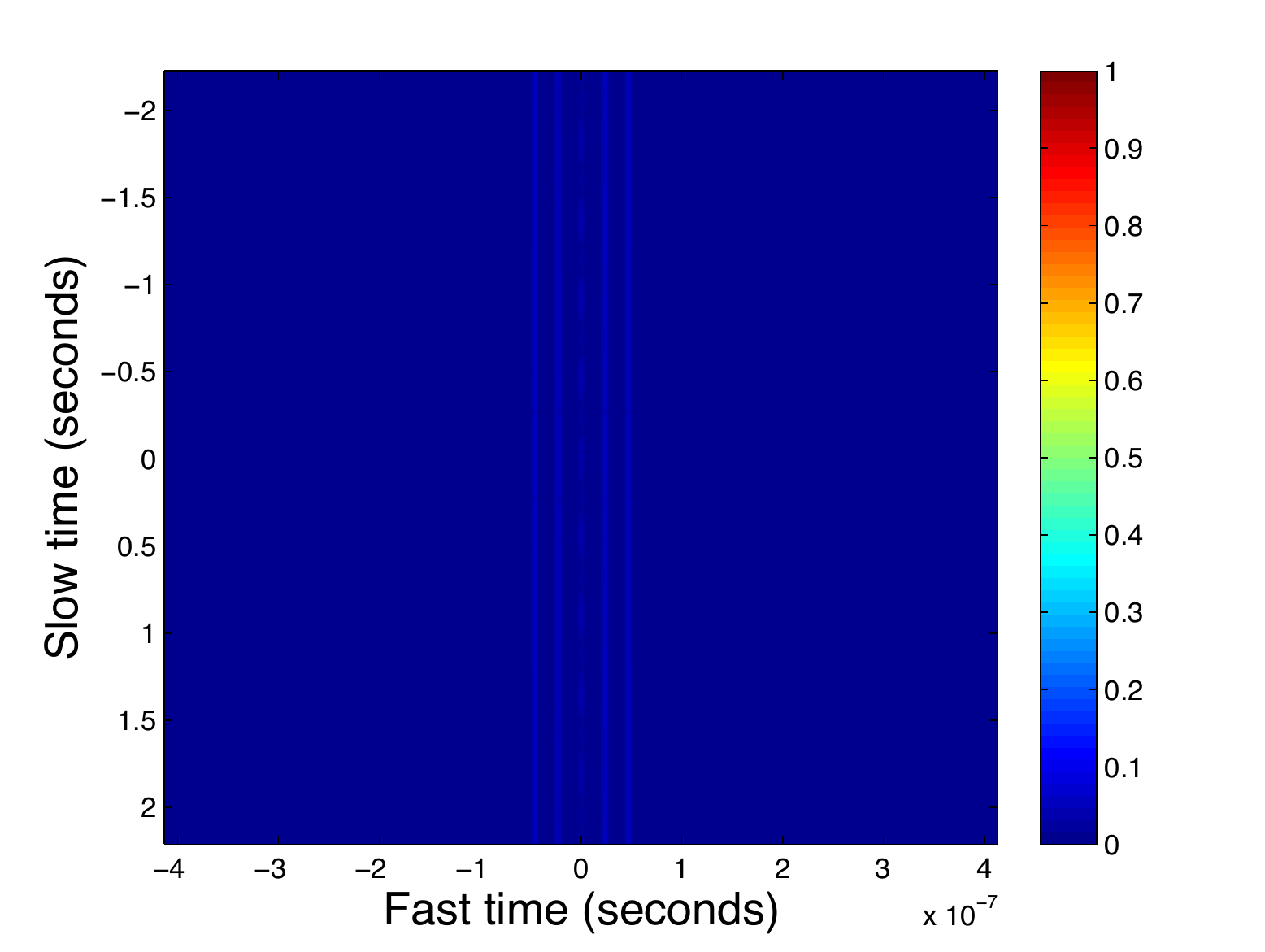}}\hspace{-0.2in}
 \subfigure{\includegraphics[width=.3\columnwidth]{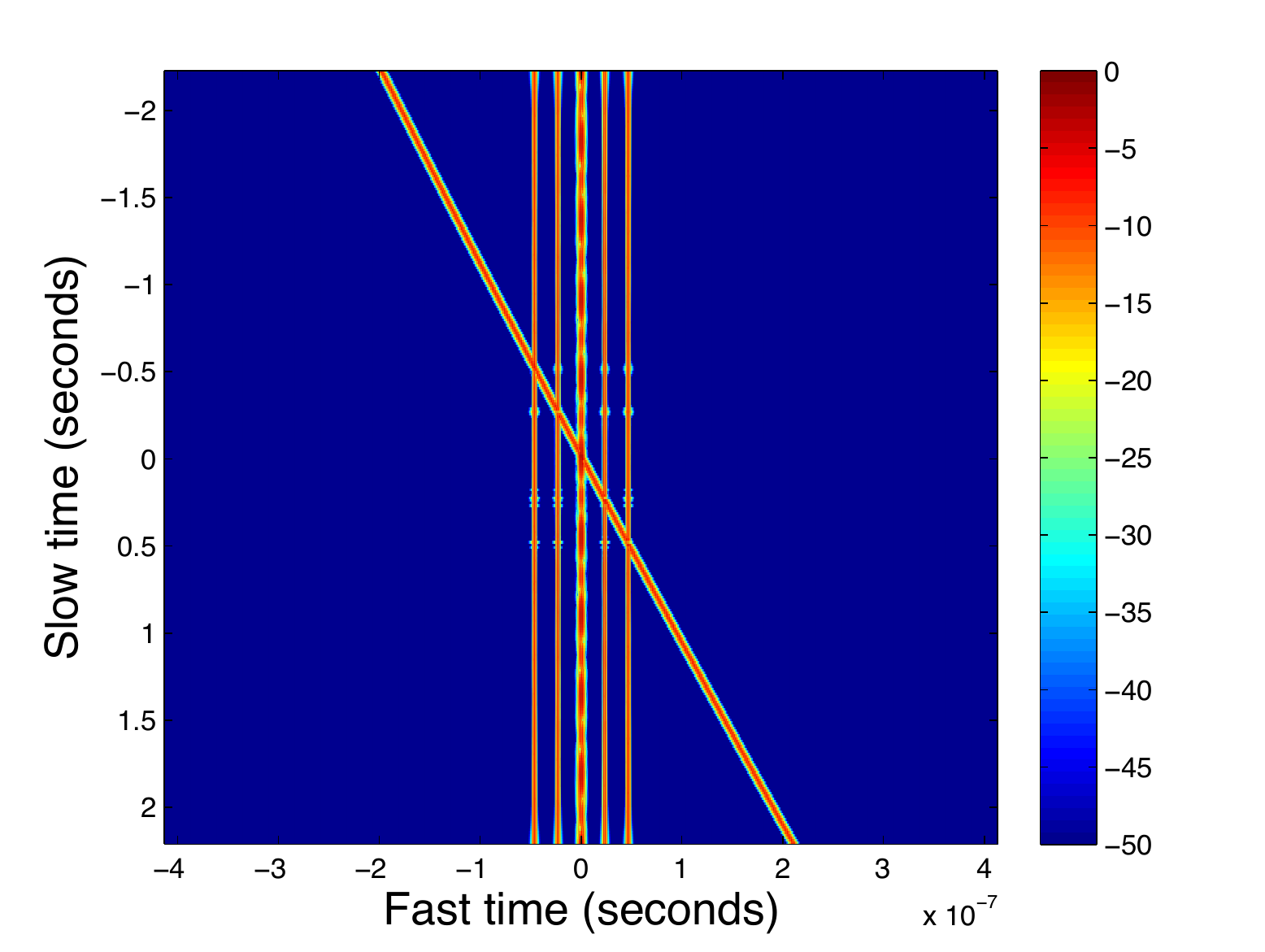}}\\
 \subfigure{\includegraphics[width=.3\columnwidth]{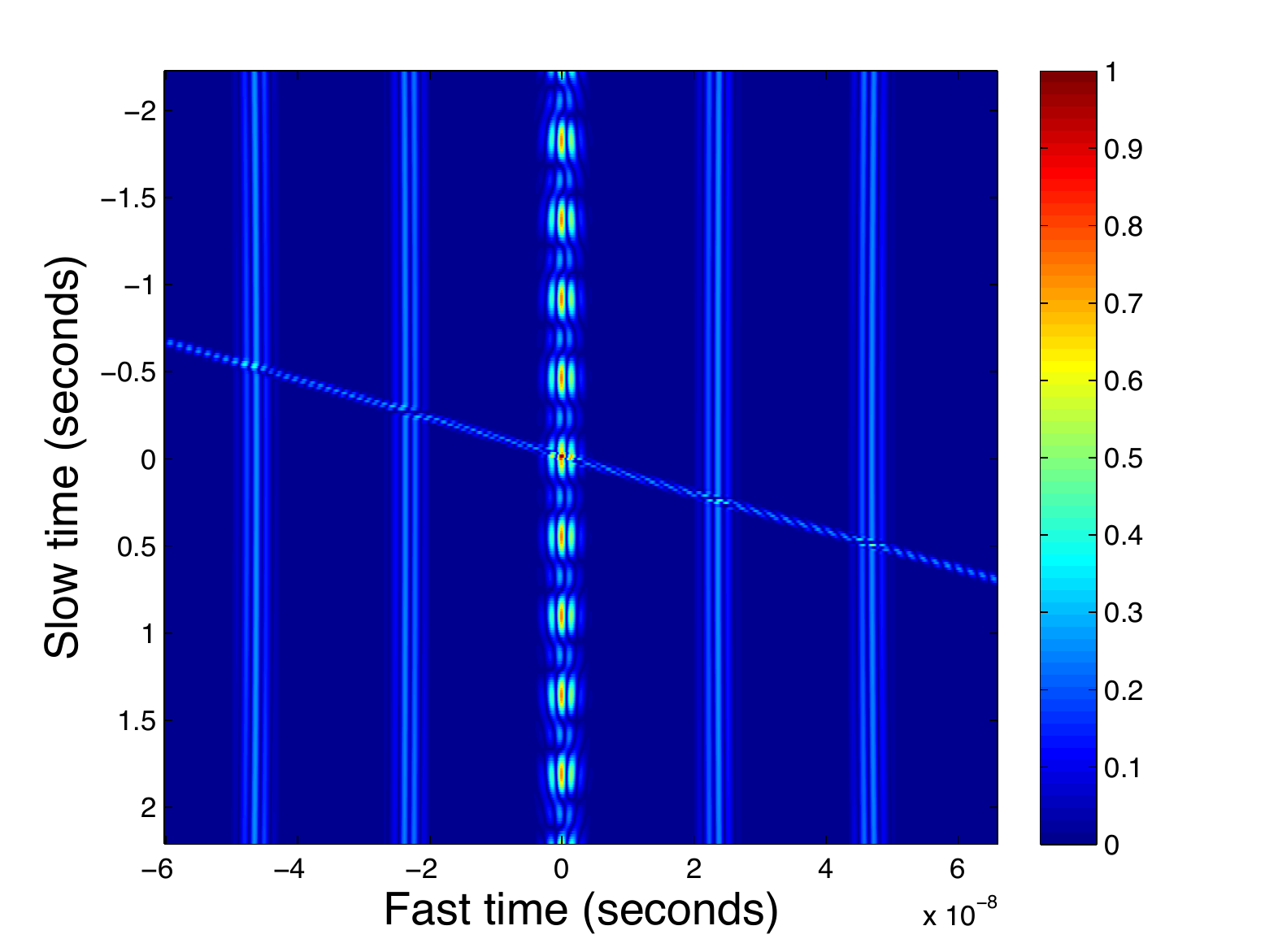}}\hspace{-0.2in}
 \subfigure{\includegraphics[width=.3\columnwidth]{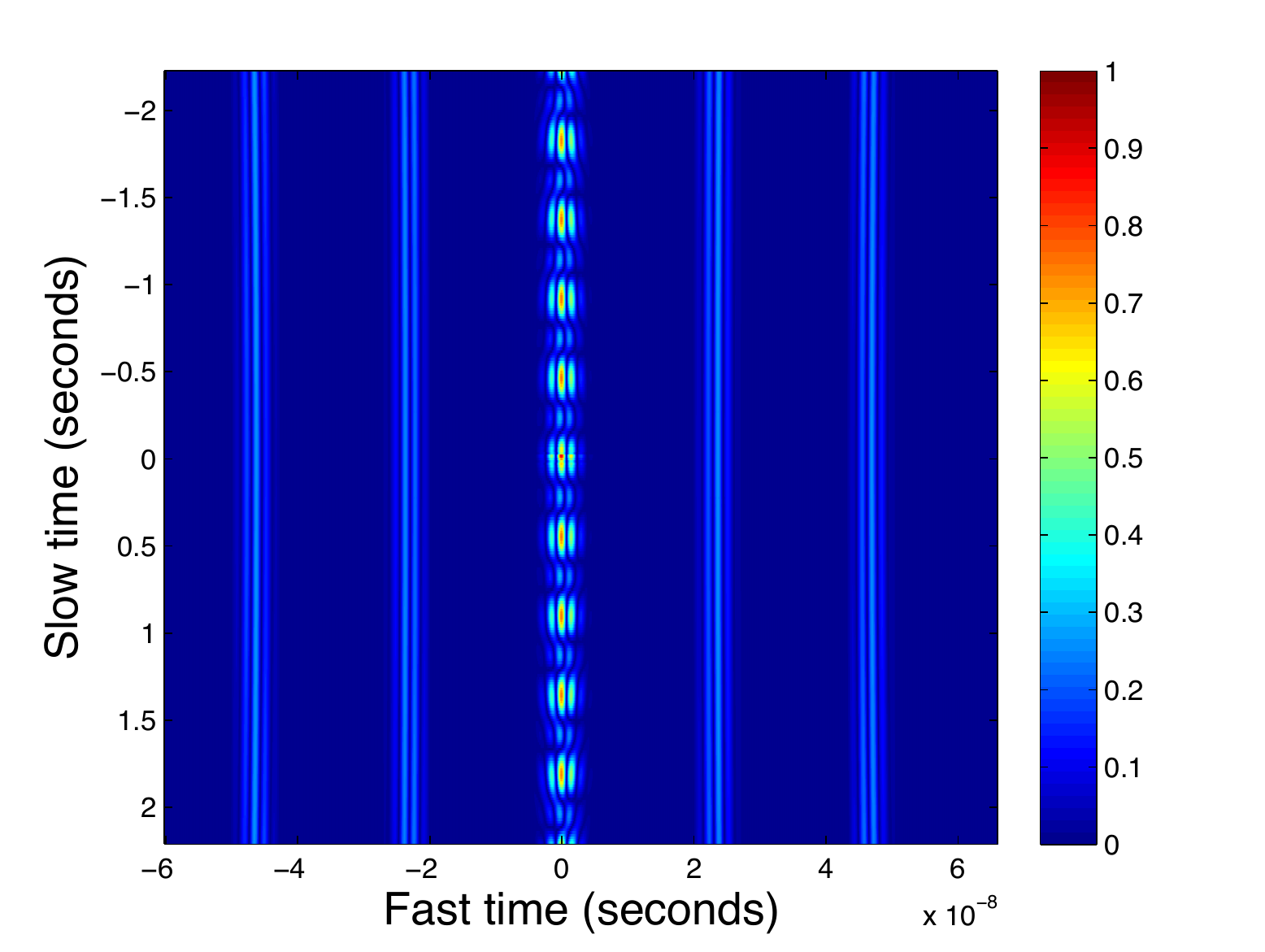}}\hspace{-0.2in}
 \subfigure{\includegraphics[width=.3\columnwidth]{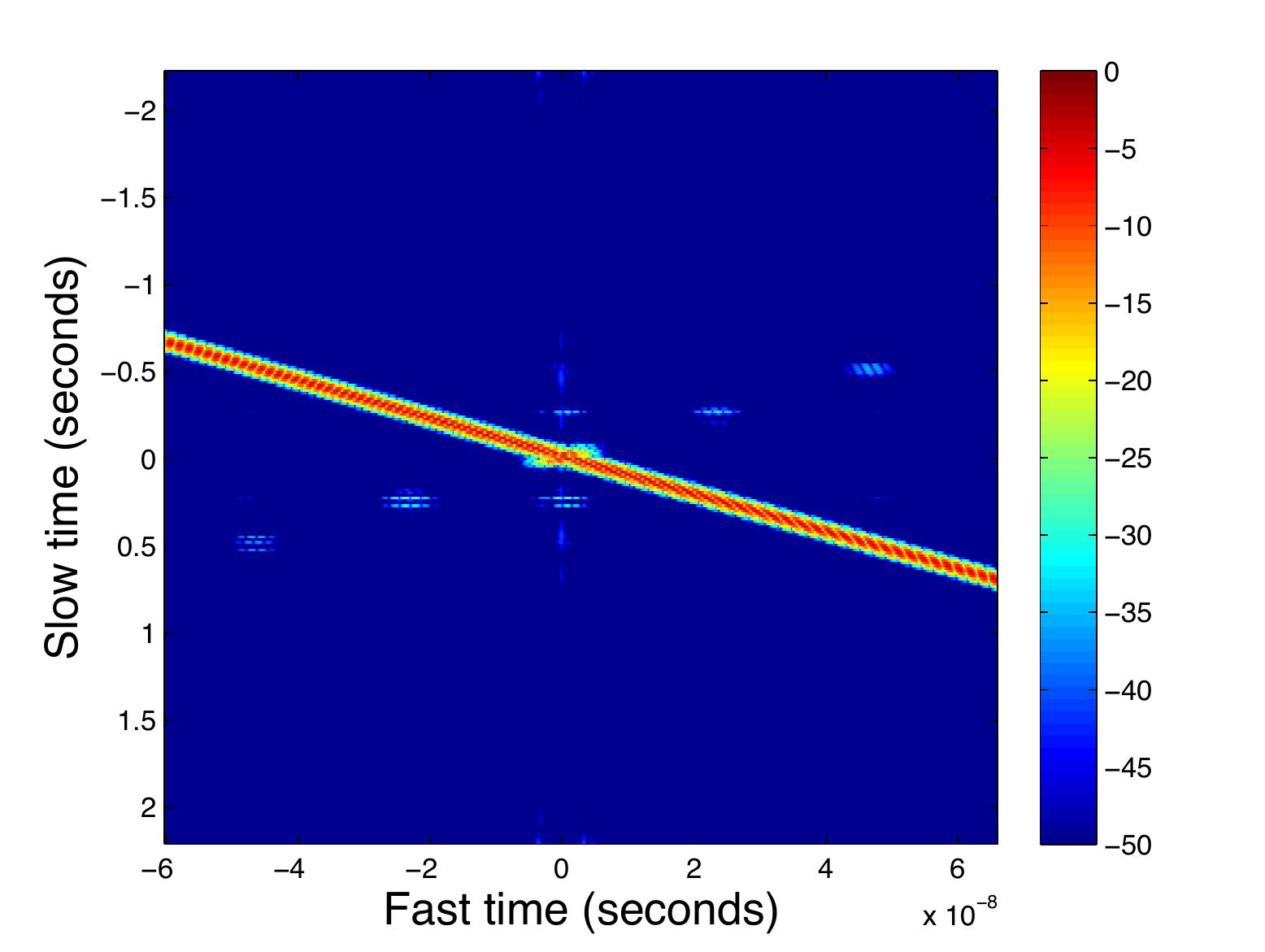}}
\end{center}
\caption{Separation of traces from the scene in Example 1.
  From left to right we show the matrix $\cM$ of traces, the low rank
  $\cL$, and sparse $\cS$ components obtained with robust PCA.  Top
  row: robust PCA applied to the entire data matrix $\cM$. Bottom row:
  robust PCA applied to a windowed $\cM$.  In each plot the matrices are
  normalized by the largest value of $|\cM|$.  The sparse components
  shown in the right column are plotted in (decibel) dB scale to make
  the contrast more visible.  }
\label{fig:rpcaEx1}
\end{figure}

\vspace{0.1in}
\noindent \textbf{Example 1:} The first example considers a scene with
one moving target and seven stationary ones.  The moving target is at
location $(0,0,0)$m at time $s=0$, corresponding to the center of the
aperture, and velocity $\vu=28/\sqrt{2}(1,1,0)$m/s.  The stationary
scatterers are located at $(0,0,0)$m, $(\pm5,0,0)$m, $(0,\pm5,0)$m,
and $(\pm10,0,0)$m. in the imaging region.  The reference point is
$\vrho_o = (0,0,0)$m.

When we apply the robust PCA method to the matrix of traces shown in
the top left plot in Figure \ref{fig:rpcaEx1}, we obtain the ``low
rank''and ``sparse'' parts shown in the top middle and top right
plots. We note from the right plot that the ``sparse'' component
contains the trace from the moving target but also remnants of the
traces from the stationary targets. The data separation is not
successful because the matrix $\cM$ is sparse to begin with.

However, if we apply robust PCA to the matrix of traces in a smaller
fast-time window, we get a successful separation.  This is shown in
the bottom row of Figure \ref{fig:rpcaEx1}.  The windowed $\cM$ has a
much larger number of non-zero entries relative to its size, and
robust PCA separates the traces from the moving target and the
stationary ones.

\begin{figure}[!t]
\centering
\includegraphics[width=.3\columnwidth]{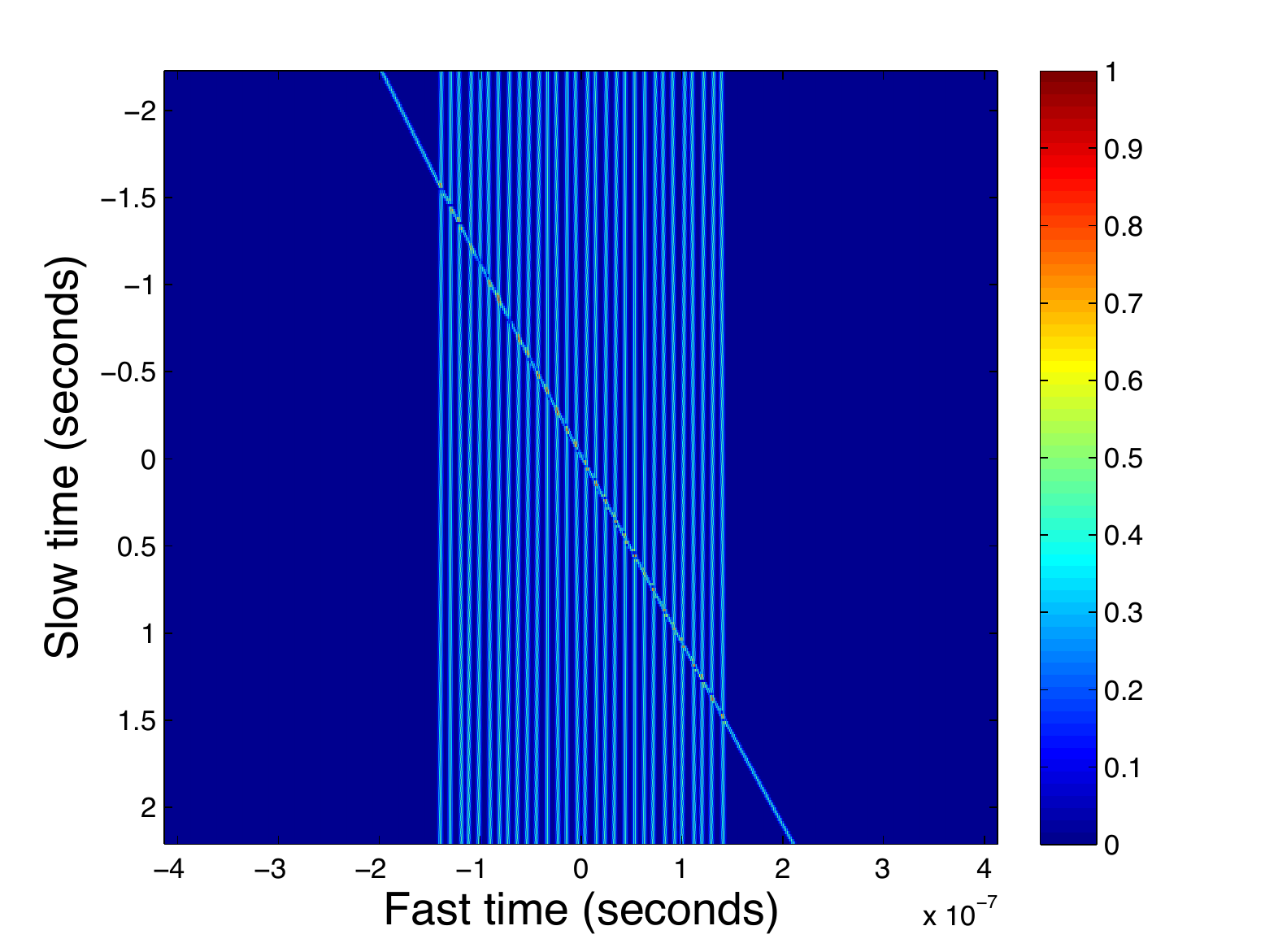}
\caption{Matrix $\cM$ of traces from the scene in Example 2.}
\label{fig:rpcaEx2a}
\end{figure}

\begin{figure}[!t]
\centering
\includegraphics[width=.9\columnwidth]{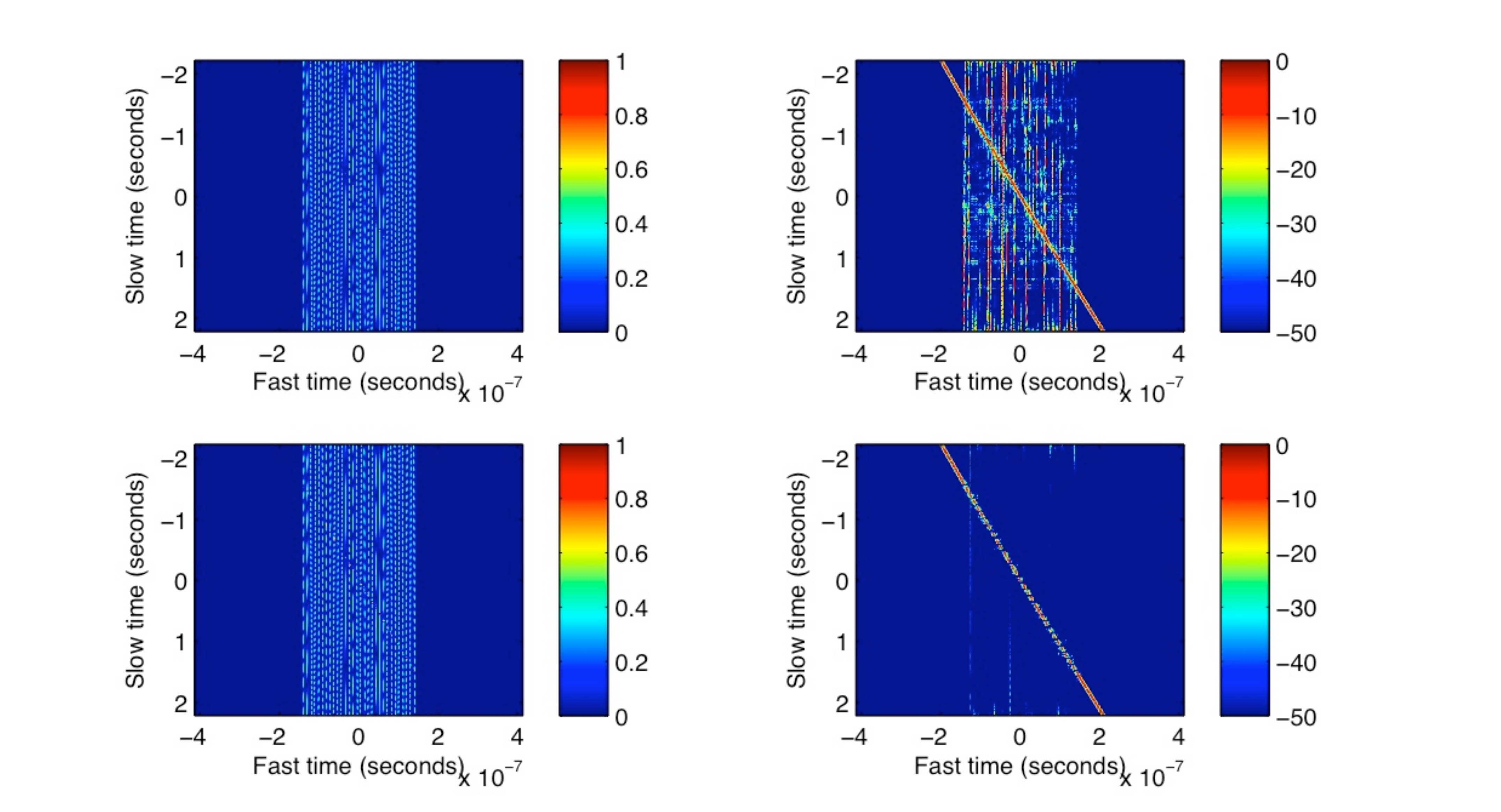}
\caption{Robust PCA separation of traces from a scene with thirty
  stationary targets and one moving target.  The top row of plots
  shows the results with robust PCA applied to the entire data matrix
  in Figure \ref{fig:rpcaEx2a}. The low rank and sparse components are
  shown on the left and right, respectively.  The bottom row of plots
  shows the results with robust PCA applied to successive small time
  windows of the data.  The sparse components shown in the right
  column are plotted in dB scale to make the contrast more visible.}
\label{fig:rpcaEx2b}
\end{figure}

\vspace{0.1in}
\noindent \textbf{Example 2:} The second example considers a scene
with thirty stationary scatterers and the same moving target as in
example 1.  The matrix $\cM$ of traces is displayed in Figure
\ref{fig:rpcaEx2a}.  The results with robust PCA applied to the entire
$\cM$ are in the top row of Figure \ref{fig:rpcaEx2b}. Again, we see
that the ``sparse'' part contains remnants of the traces from the
stationary targets. This is not because $\cM$ is sparse, as was the
case in example 1, but because the traces from the stationary targets
do not form a matrix of sufficiently low rank. When we window $\cM$,
we work with a few nearby traces at a time. These traces are similar
to each other so they form a low rank matrix that can be separated
with robust PCA. The final result is given by the concatenation of the
matrices in each window. It is shown in the bottom row of Figure
\ref{fig:rpcaEx2b}.

These two examples show that windowing is a key step to successful
data separation with robust PCA.  To choose appropriate window sizes,
it is necessary to understand how the rank and sparsity of the data
matrix $\cM$ depend on the location and density of the stationary
targets, and on the velocity of moving targets. This is addressed in
the analysis of \cite{BCP12}, which we review in section
\ref{sect:rpcaanalysis}.

% ------------------------------------------
\subsection{Discussion}
\label{sect:sepdisc}
\begin{figure}[!h]
\centering
\includegraphics[width=.8\columnwidth]{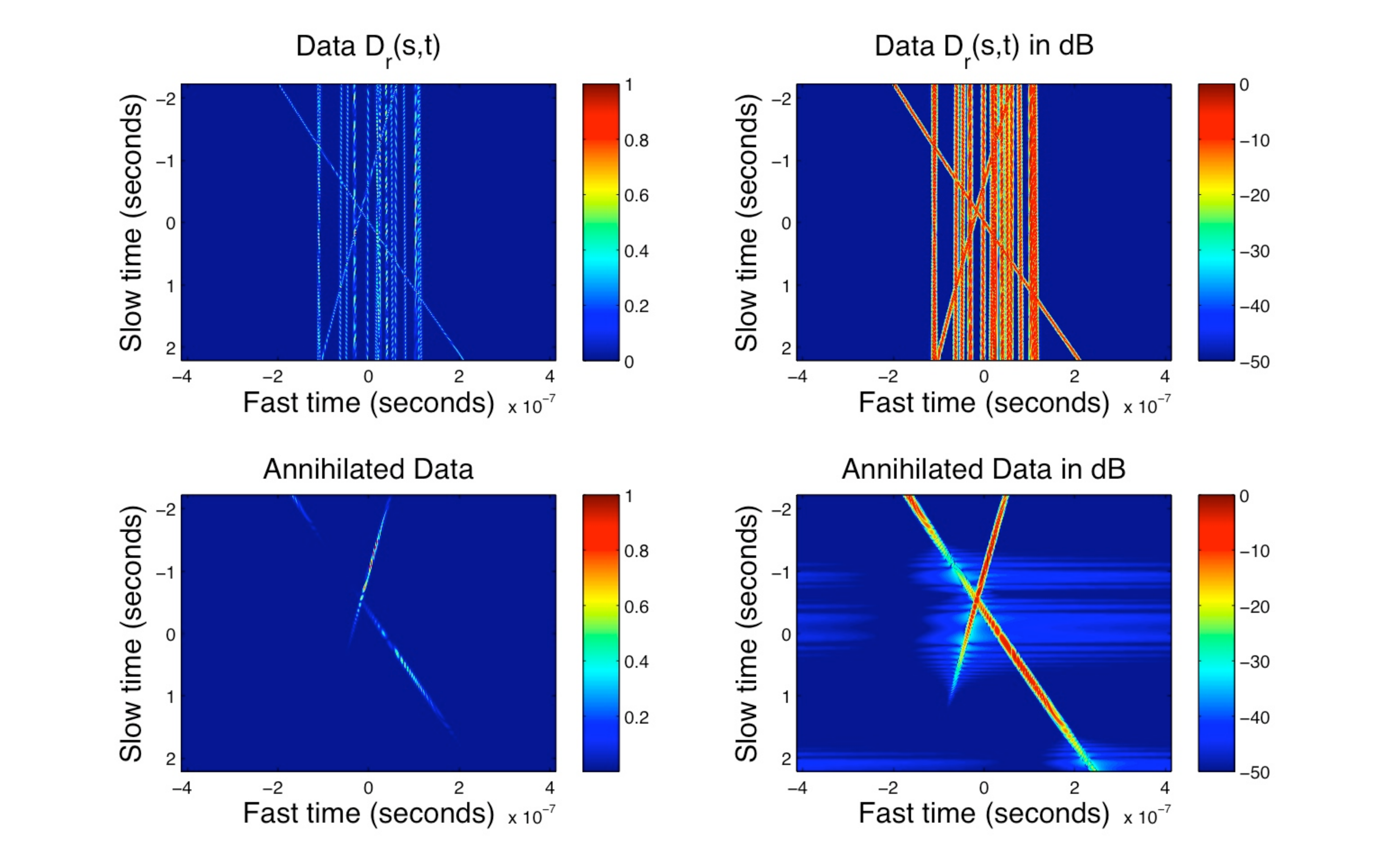}
\caption{Annihilation with the filter (\ref{eq:OP2}) of the stationary
target traces from Scene 1.  Top: the matrix of traces. Bottom: the
filtered traces. In each plot we show absolute values normalized by
the largest value of $|\cM|$. The right column is the left column
plotted in dB scale.}
\label{fig:annihilScene1}
\end{figure}

\begin{figure}[!t]
\centering
\subfigure{\includegraphics[width=0.3\columnwidth]{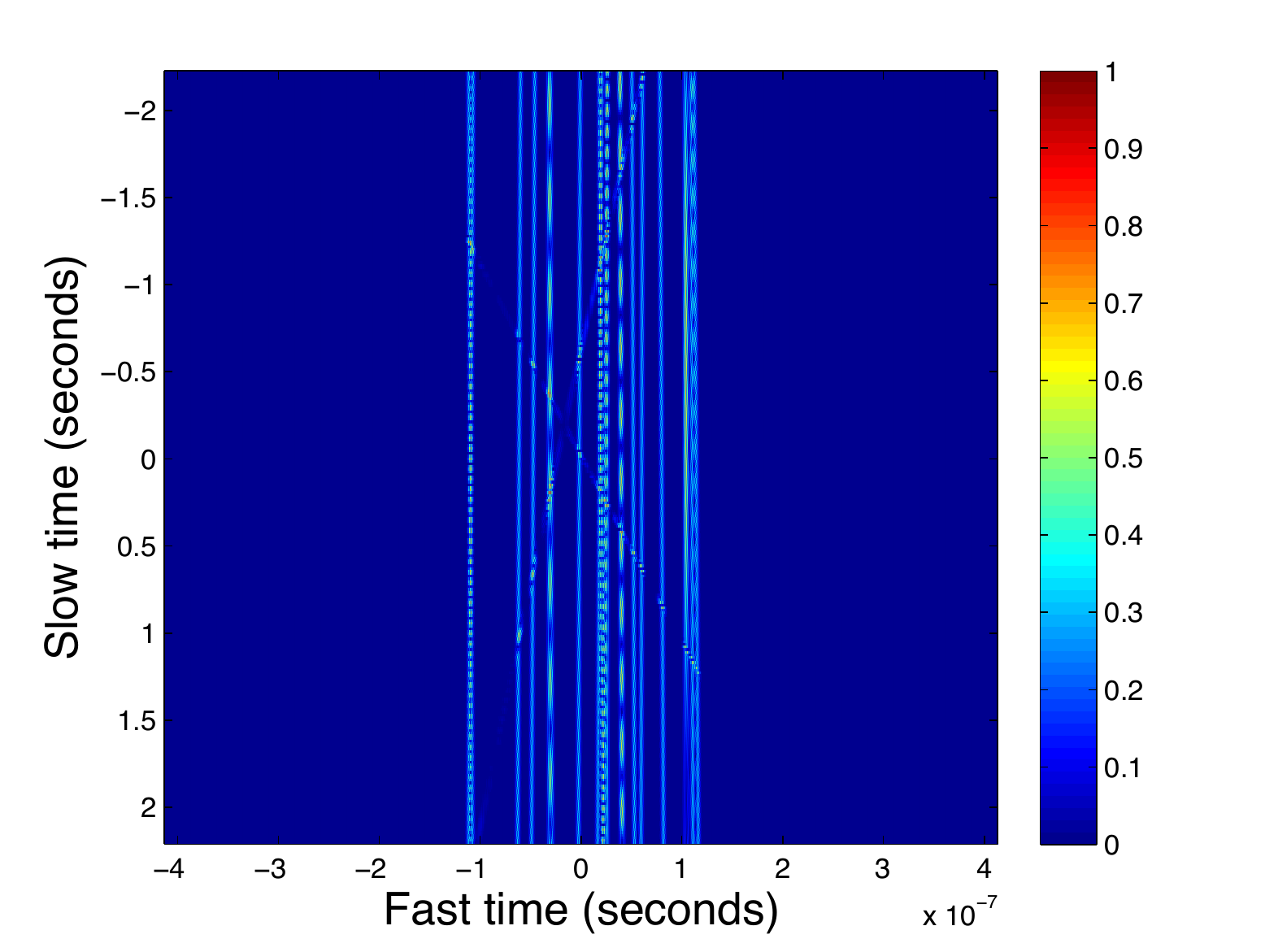}}
\subfigure{\includegraphics[width=0.3\columnwidth]{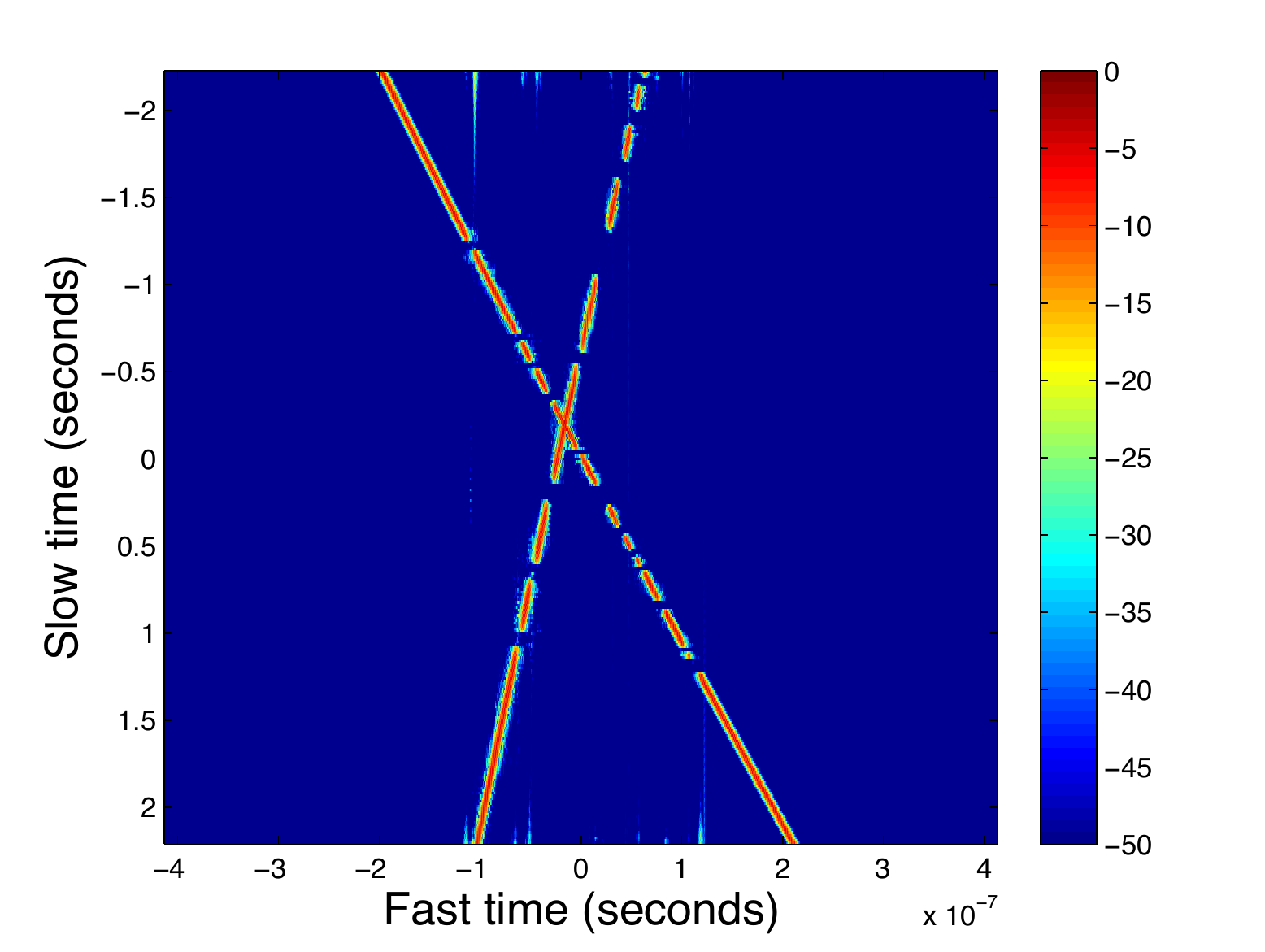}}
\caption{Trace separation with robust PCA for Scene 1. The low rank
part $\cL$ is shown on the left and the sparse part $\cS$ is shown on
the right.  In each plot we show absolute values normalized by the
largest value of $|\cM|$.  The sparse part is plotted in dB scale.}
\label{fig:rpca2mov}
\end{figure}

We use numerical simulations for Scene 1 to compare the two
methods described in sections \ref{sect:ANNIH} and \ref{sect:RPCA}.
The results are in Figures \ref{fig:annihilScene1} and
\ref{fig:rpca2mov}.


The first method uses the filter $\mathcal Q$ defined in
(\ref{eq:OP2}) to annihilate the traces from the stationary
targets. We see from Figure \ref{fig:annihilScene1} that this is
accomplished, but the filter removes some parts of the traces from the
moving targets as well. All our numerical simulations show that
$\mathcal Q$ separates nicely the traces from the moving targets if the
scene does not have too many stationary targets at almost the same
range. Otherwise, it also removes parts of the traces from the
moving targets.

The data separation with the robust PCA approach is better, as seen in
Figure \ref{fig:rpca2mov}. This method is more sophisticated than the
other, because it involves careful windowing of the traces, but it
works better for complex scenes that are not too large.  Robust PCA
relies heavily on range compression removing most of the slow-time
dependence of the traces from the stationary targets. But for extended
scenes there is no single reference point $\vrho_o$ that accomplishes
this task, and data separation with robust PCA deteriorates. The
filter (\ref{eq:OP2}) still works for extended scenes, because the
travel time transformations in (\ref{eq:OP1}) are relative to the
location of the targets that we wish to annihilate, and not some
arbitrary reference point.

In principle, we can use a combination of the two approaches to
improve the data separation with robust PCA for extended scenes. The
idea is to begin with a preliminary image and identify roughly the
stationary targets to be separated from the scene.  Then we can window the
image around groups of such targets and work with one subscene at a
time. It is important to observe here that if the aperture is not too
large, as is the case in motion estimation, the image is approximately
the Fourier transform of the data with a phase correction. See for
example \cite{sar} for a detailed explanation of this fact. Thus, we
can obtain approximately the traces from the smaller subscene by
Fourier transforming the windowed image. To separate these traces, we
can apply the travel time transformation (\ref{eq:TT+}) with a
reference point $\vrho^e$ in the subscene, followed by the robust PCA
approach as described above.

% ------------------------------------------
\subsection{Velocity estimation and separation of traces from 
multiple moving targets}
\label{sect:rotationopt}
Let us assume in this section that the traces from the stationary
targets have been removed, and denote by the same $D_r(s,t)$ the
traces from the moving targets.  Since there may be multiple moving
targets, we would like to separate these traces further, based on
differences between the target velocities. We can do so by extending
the annihilation filters described in section \ref{sect:ANNIH}, as we
explain next.

We need a key observation, shown with analysis in section
\ref{sect:analysis}. The slope of the trace from a moving target is
determined approximately by its speed along the unit
vector
\begin{equation}
\vec{\bf m}_o = \frac{\vr(0)-\vrho_o}{|\vr(0)-\vrho_o|},
\label{eq:Mo}
\end{equation}
pointing from the SAR platform location $\vr(0)$ at the middle of the
aperture to the reference point $\vrho_o$ in the imaging
scene. Moreover, the curvature of the trace is determined by the speed
projected on the plane orthogonal to $\vec{\bf m}_o$.

Because we assume a flat imaging surface, the target speed $\vec{\bf
  u}$ is a vector in the horizontal plane
\begin{equation}
\vu = (\bu,0), \quad \bu \in \mathbb{R}^2.
\end{equation}
We parametrize it by the range velocity $u$ and the cross-range
velocity $u_\perp$.  The range velocity is 
\begin{equation}
u = \vec{\bf u} \cdot \vec{\bf m}_o = \bu \cdot \bm_o,
\label{eq:rangeu}
\end{equation}
where $\bm_o$ is the projection of $\vm_o$ on the imaging plane.
The cross-range velocity is defined by 
\begin{equation}
u_\perp = \vec{\bf t} \cdot \mathbb{P}_o \vec{\bf u} = \bu \cdot 
\bt - u \, \vm_o \cdot \vt(0),
\label{eq:crangeu}
\end{equation}
where $\mathbb{P}_o$ is the orthogonal projection 
\begin{equation}
\label{eq:ORTHPROJ}
\mathbb{P}_o = I - \vm_o \vm_o^T,
\end{equation}
and $\vec{\bf t}$ is the unit vector tangent to the SAR platform
trajectory, at the center $\vr(0)$ of the aperture. Note that in many
setups $\vm_o$ and $\vt$ are almost orthogonal, so $u_\perp$ is
approximately the target speed along the vector $\bt$, the
projection of $\vt$ on the imaging plane.

Next, we define the transformation mapping $\mathbb{T}^{\vu^e,
\vrho^e}_+$, the generalization of (\ref{eq:TT+}),
\begin{equation}
\left[\mathbb{T}^{\vu^e, \vrho^e}_+ D_r \right](s,t) = D_r\left(s,t +
\Delta \tau(s,\vrho^e + s \vu^e)\right),
\label{eq:TTu+}
\end{equation}
with $\Delta \tau$ given in (\ref{eq:DeltaTau}). Here $\vrho^e$ is the
estimated location of the moving target at the origin of the slow
time, and $\vu^e$ is the estimated velocity. If the line segment
\[
\vrho^e + s \vu^e, \quad 
s \in [-S(a),S(a)],
\]
is close to the target trajectory, then the transformation
(\ref{eq:TTu+}) removes approximately the dependence of its trace on
the slow time.  Thus, we can subtract it from the other traces using
the same difference operator $\mathbb{D}_s$.  The annihilation filter
is the generalization of (\ref{eq:OP1}),
\begin{eqnarray}
\left[ \mathcal Q^{\vu^e, \vrho^e} D_r \right](s,t) &=& \left[
\mathbb{T}^{\vu^e, \vrho^e}_- \mathbb{D}_s \mathbb{T}^{\vu^e,
\vrho^e}_+ D_r \right] (s,t) \nonumber \\ &=& \left[\mathbb{D}_s
D_r\left(s,t'+ \Delta \tau(s,\vrho^e + s \vu^e)\right) \right]_{t'= t-
\Delta \tau(s,\vrho^e + s \vu^e)},
\label{eq:OPu1}
\end{eqnarray}
where the mapping $\mathbb{T}^{\vu^e, \vrho^e}_-$ undoes the travel
time transformation.  Additionally we can apply robust PCA to 
$\left[\mathbb{T}^{\vu^e, \vrho^e}_+ D_r \right](s,t)$.  This should result
in the low rank matrix containing the trace corresponding to the moving
target with velocity $\vrho^e$ and the sparse matrix containing the rest
 of the moving target traces.

The implementation of (\ref{eq:TTu+}) and (\ref{eq:OPu1}) requires the estimates $\vrho^e$
and $\vu^e$.  There are two scenarios. The first assumes that we know
$\vrho^e$ from tracking the target in previous sub-apertures, and
seeks only $\vu^e$. The second seeks both $\vrho^e$ and $\vu^e$. In
either case, the estimation begins with that of the range velocity,
which controls the slope of the trace. Then, we can estimate $\vrho^e$
if we need to, and the cross-range speed $u_\perp^e$.

We estimate $u^e$ as the maximizer of the objective function
\begin{equation}
\label{eq:OBJg}
g(u) = \max_{-m/2 \le \ell \le m/2} \sum_{j=-n/2}^{n/2}
\left|\left[\mathbb{T}^{\vu,\vrho^e}_+D_r\right](s_j,t_\ell)
\right|, \quad \vu  = u \, (\bm_o,0).
\end{equation}
 The idea is that when $u^e$ is close to the range speed of the
target, the transformation $\mathbb{T}^{\vu,\vrho^e}_+$ makes the
trace an approximate line segment, and when we sum it over $s$, we get
a peak around some fast time $t_\ell$, for $|\ell| \le m/2$.  If we do not
know $\vrho^e$, we can set it equal to $\vrho_o$ in
(\ref{eq:OBJg}). Then, we can form an image $\cI_{\vu^e}(\vrhoi)$
using (\ref{eq:imagecompu}) with $\vu^e = u^e (\bm_o,0)$, and estimate
$\vrho^e$ from there. The motion compensation in (\ref{eq:imagecompu})
brings the moving target approximately in focus, and this is why we
expect to be able to estimate its location from the image.

Once we have estimated $u^e$ and $\vrho^e$, we seek the cross-range
speed $u_\perp^e$ as the minimizer of the objective function 
\begin{equation}
g_\perp(u_\perp;u^e, \vrho^e)=\sum_{\ell=-m/2}^{m/2}\sum_{j=-n/2}^{n/2}
\left|\mathbb{D}^2_s \left[\mathbb{T}^{\vu,\vrho^e}_+
D_r\right](s_j,t_\ell)\right|,
\label{eq:OBJ2}
\end{equation}
which quantifies the curvature of the trace.  Here $\mathbb{D}^2_s$
approximates the second derivative in $s$, and $\vu = (\bu,0)$ is
defined by
\begin{equation}
\bu \cdot \bm_o = u^e, \quad \bu \cdot \bt(0) = u_\perp + u^e \vm_o
\cdot \vt(0).
\end{equation}
Because $u_\perp$ has a weaker effect on the traces than $u$, it is
more difficult to estimate it. Even if there is no noise, the
estimation is sensitive to remnants of the traces from the stationary
targets that have not been perfectly annihilated, the traces from the
other moving targets, and the discrete sampling of the data. If we
took the maximum over $\ell$ of the sum
\[
\cF(u_\perp, \ell) = \sum_{j=-n/2}^{n/2}
\left|\left[\mathbb{D}^2_s \mathbb{T}^{\vu,\vrho^e}_+
D_r\right](s_j,t_\ell)\right|
\]
as we did in the objective function (\ref{eq:OBJg}), we would seek a
saddle point of $\cF$ to estimate $u_\perp$. That is to say, we would
minimize with respect to $u_\perp$ and maximize with respect to
$\ell$. Our numerical simulations show that such an estimation is not
robust.  To stabilize it, we sum over $\ell$ in (\ref{eq:OBJ2}) and
reduce the problem to that of minimizing in $u_\perp$.

\begin{figure}[!t]
\centering
\includegraphics[width=1.\textwidth]{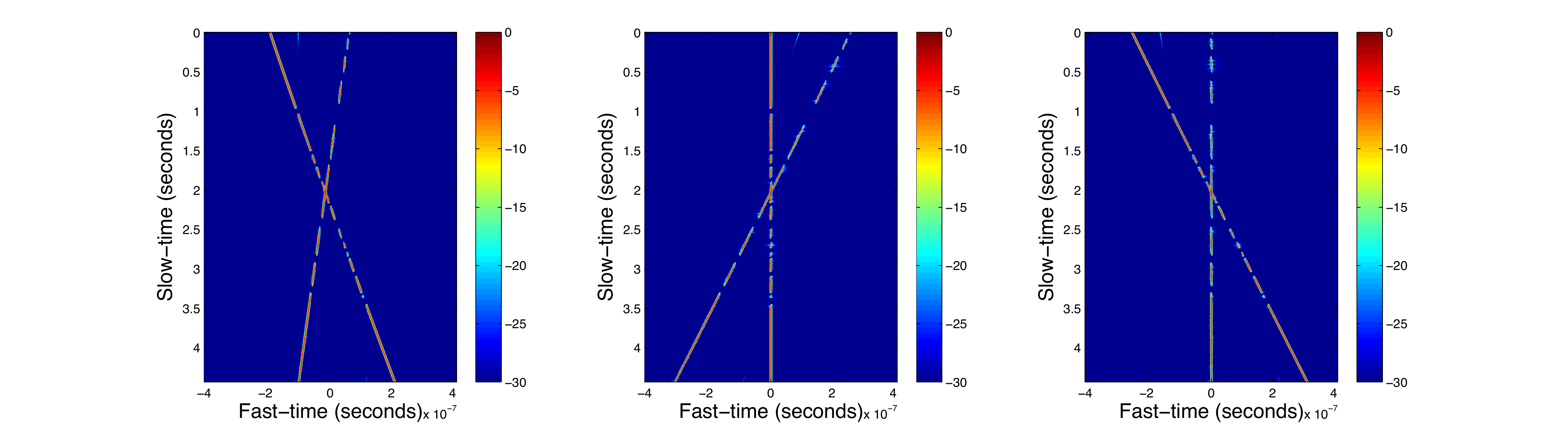}
\caption{Left: The traces from the two moving targets in Scene 1.
Middle and right: The transformed traces $\mathbb{T}^{\vu,\vrho}_+
D_r(s,t)$ for the first and second target.  }
\label{fig:rotateTrueFilt}
\end{figure}
\begin{figure}[!t]
\centering
{\includegraphics[width=.6\textwidth]{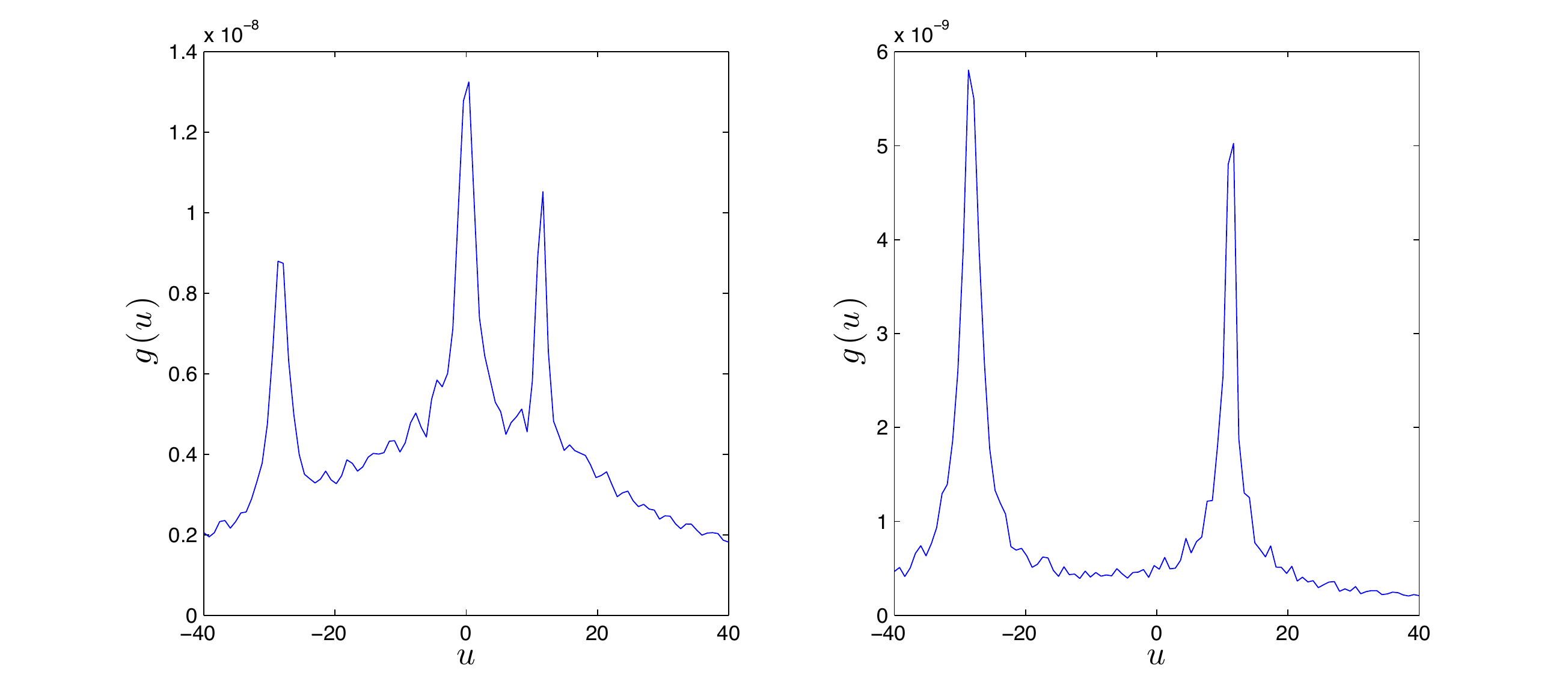}}
\caption{The objective function $g(u)$ using the traces from the
entire Scene 1 (left) and the traces from the moving targets in Scene
1 (right).}
\label{fig:filtComp}
\end{figure}

As an illustration, we present the results for Scene 1. Figure \ref{fig:rotateTrueFilt}
displays $\mathbb{T}^{\vu,\vrho}_+D_r(s,t)$ using the true locations 
$\vrho$ and velocities $\vu$ for the two moving targets. Note how after
the transformation the traces become vertical line segments. 
If we were to search for the range speed $u^e$ without annihilating the traces of
the stationary targets, we would get the objective function $g(u)$
plotted on the left in Figure \ref{fig:filtComp}. The prominent peak of
$g$ is at speed $u = 0$, because the stationary targets dominate the
scene. The smaller peaks are due to the two moving targets.  When we
compute the objective function $g(u)$ with the traces from the moving
targets alone, we get two large peaks, at the range speeds of the targets. 
%The effect of the travel time transformation (\ref{eq:TTu+}) on the traces is illustrated in Figure \ref{fig:rotateTrueFilt}.
 The speed in (\ref{eq:TTu+}) is $\vu = u^e (\bm_o,0)$, with $u^e$ estimated
from the first and second peak of $g(u)$, respectively. 

\begin{figure}[!h]
\centering
{\includegraphics[width=.7\textwidth]{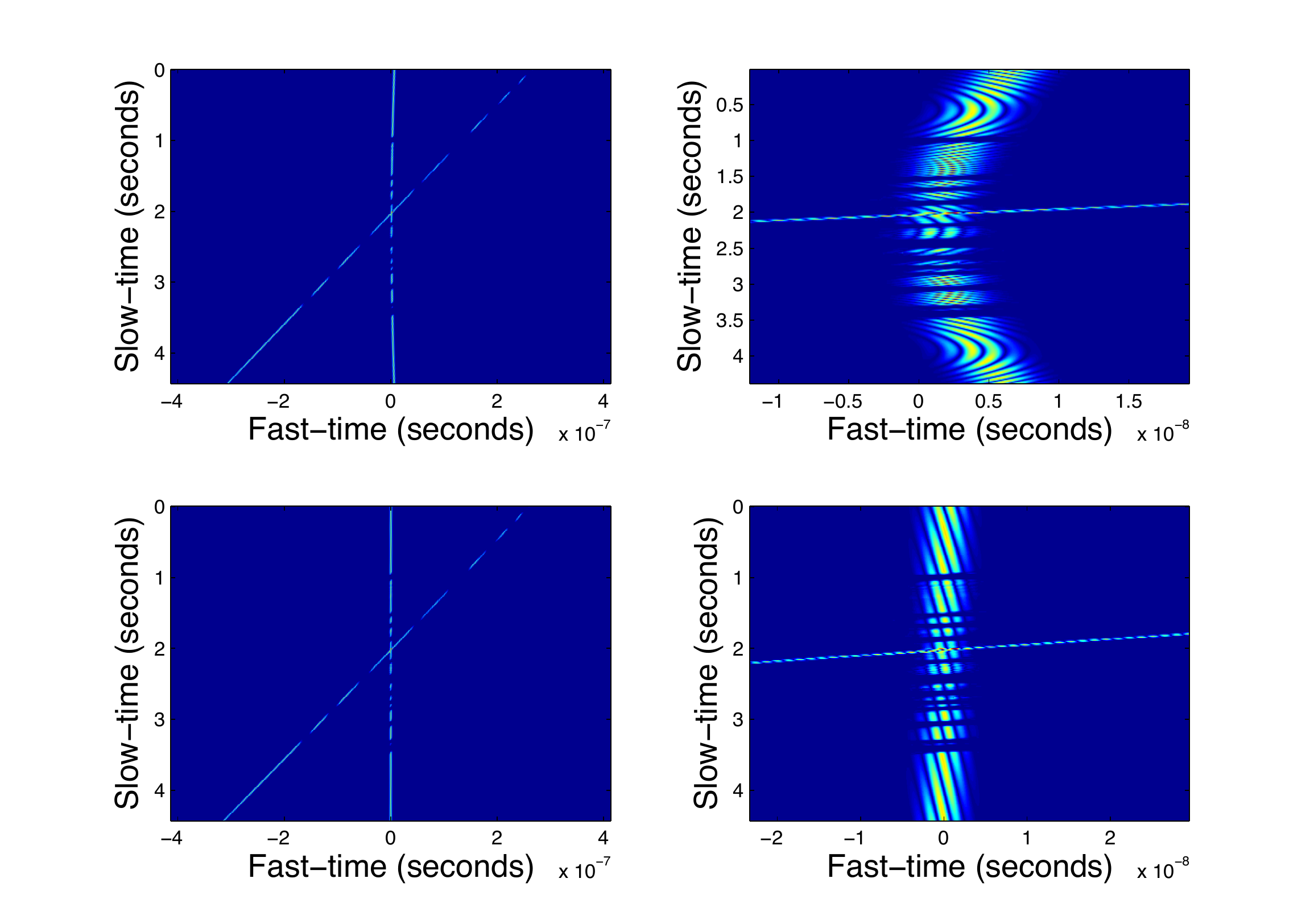}}
\caption{Top row: The traces after the travel time transformation
(\ref{eq:TTu+}) at the estimated range speed of one target. 
Bottom row: The traces after the travel time transformation
(\ref{eq:TTu+}) with the estimated cross-range correction of the speed. 
The right column is a zoom of the left column.} 
\label{fig:rotate2stage}
\end{figure}

\begin{figure}[!h]
\centering
{\includegraphics[width=.7\textwidth]{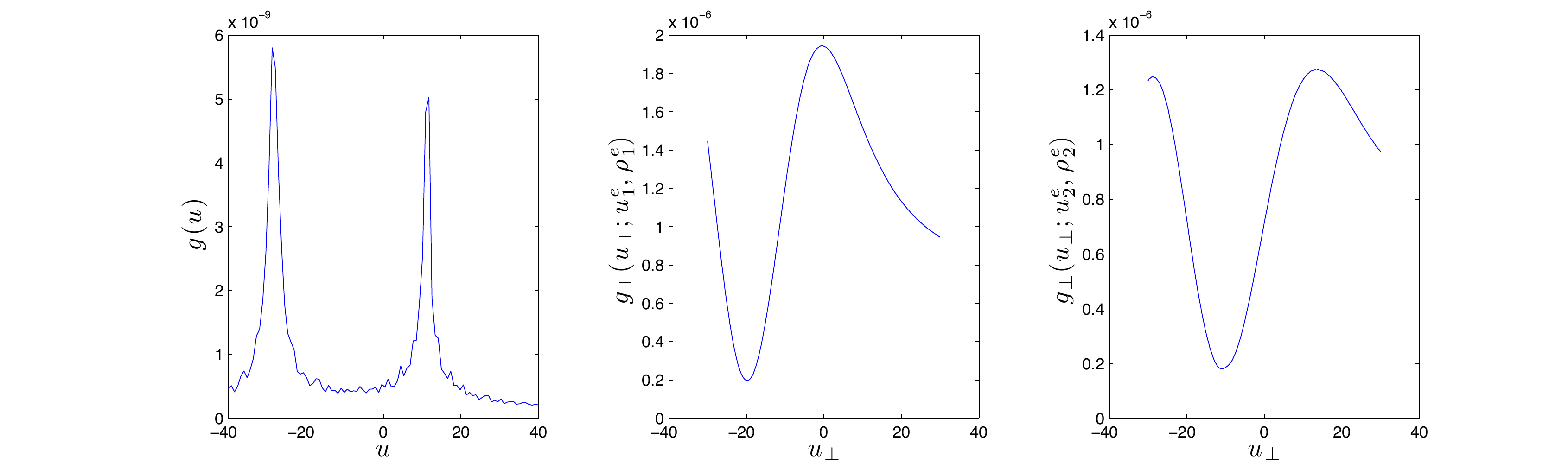}}
\caption{Left: The objective function $g$.  Middle and right: The
  objective function $g_\perp$ for the estimated range speeds of the
  two targets in Scene 1. The minima of $g_\perp$ are the estimates
  of $u_\perp$.}
\label{fig:objfxns}
\end{figure}

Recall that the cross-range speed $u_\perp$ has a weaker effect on the
traces than the range speed, it determines only their curvature.  We
illustrate this in the top right plot of Figure
\ref{fig:rotate2stage}, where we zoom around the trace of one target,
rotated with transformation (\ref{eq:TTu+}) at the estimated range
speed $u^e(\bm_o,0)$.  The bottom row in Figure \ref{fig:rotate2stage}
shows the traces transformed with (\ref{eq:TTu+}), after the
cross-range speed correction. The trace is now straightened, as seen
in the bottom right plot.  The objective function (\ref{eq:OBJ2}) is
shown in the middle and right plots of Figure \ref{fig:objfxns}. We
plot it for the two range speeds, the peaks of $g$ shown on the
left. The function $g_\perp$ has minima at the cross-range speeds of
the two targets.

% ------------------------------------------
\section{Analysis}
\label{sect:analysis}
\setcounter{equation}{0} In this section we give a brief analysis of
our data separation approaches. The analysis is for simple scenes, but
the methods work in complex scenes, as shown already in the previous
section. More numerical results are in section \ref{sect:numeric}.  We
begin with the data model and the scaling regime. Then, we analyze the
annihilation filters defined in section \ref{sect:ANNIH}.  The
analysis extends to the filters described in section
\ref{sect:rotationopt}, once we explain how the target velocity
affects the traces. We end the section with a brief review of the
analysis in \cite{BCP12} of the data separation with robust PCA.

% ------------------------------------------
\subsection{Scaling Regime}
\label{sect:scales}
We illustrate our scaling regime using the GOTCHA Volumetric SAR data
set.  The relevant scales are the central frequency $\nu_o$, the
bandwidth $B$, the typical range $L$ from the SAR platform to the
targets, the aperture $a$, the magnitude $|\vu|$ of the velocity of
the targets, the speed $V$ of the SAR platform and $R^{\cI}$, the
diameter of the imaging set $\cJ^{^\cI}$. We let $L=|\vr(0) -
\vrho_o|$, and assume that 
\begin{equation}
B \ll \nu_o, \quad a \ll L, \quad R^{^\cI} \ll L.
\label{eq:SC1}
\end{equation}
We also suppose that the target speed is smaller than that of the SAR
platform
\begin{equation}
|\bu| < V.
\label{eq:SC2}
\end{equation}
The GOTCHA parameters satisfy these assumptions:
\begin{itemize}
\item The central frequency is $\nu_0=9.6$GHz and the bandwidth is
$B=622$MHz.
\item The SAR platform trajectory is circular, at height $H=7.3$ km,
with radius $R = 7.1$km and thus $L\approx 10$km. We consider
sub-apertures of one circular degree, which corresponds to $a=124$m,
and image in domains of radius $R^{^\cI} \le 50$m.
\item The platform speed is $V=70$m/s and the targets move with
velocities $|\vu|\le 28$m/s.
\item The pulse repetition rate is $117$ per degree, which means that
a pulse is sent every $1.05$m, and $\Delta s = 0.015$s.
\end{itemize}

For a stationary target, we obtain from basic resolution theory
\cite{sar,Jakowatz} that the range resolution is $c/B=48$cm, and the
cross range resolution is $\lambda_0 L/a = 2.5$m, with one degree
aperture $a$ and central wavelength $\la_0 = 3$cm.

% ------------------------------------------
\subsection{Data model}
\label{sect:model}
We use a simple data model that approximates the scatterers by point
targets in the horizontal plane and neglects multiple scattering
between them. We also suppose that the target trajectories can be
approximated locally, for the slow time interval defining a single
sub-aperture, by the straight line segment
\[
\vrho(s) = \vrho(0) + s \vu,
\]
with constant speed $\vu = (\bu,0)$, and $\bu \in \mathbb{R}^2$.  The
traces are modeled by
\begin{equation}
  D_r(s,t)\approx \sum_{q = 1}^{N}\frac{\sigma_q(\om_o)}{
    (4\pi|\vr(s)-\vrho_q(s)|)^2}f_p(t-(\tau(s,\vrho_q(s))-
    \tau(s,\vrho_o))),
\label{eq:MOD1}
\end{equation}
where $N$ is the number of targets. They are located at $\vrho_q(s)$
and with reflectivity $\sigma_q(\om_o)$, for $q = 1, \ldots, N$.  The
model neglects the displacement of the targets during the round trip
travel time. This is justified in radar applications because the waves
travel at the speed of light which is several orders of magnitude
larger than the speed of the targets. We refer to \cite{sar} for a
derivation of model (\ref{eq:MOD1}) in the scaling regime described
above. We simplify it further by approximating the amplitude factors
\[
\frac{1}{4 \pi |\vr(s)-\vrho_q(s)|} \approx \frac{1}{4 \pi L},
\]
and assuming that the targets have the same reflectivity
\[
\sigma_q(\omega_o)=\sigma(\omega_o),\quad q=1,\ldots,N.
\]

The model of the traces becomes
\begin{equation}
D_r(s,t)\approx \frac{\sigma(\omega_o)}{(4\pi L)^2}\mathcal M(s,t),
\label{eq:model}
\end{equation}
with
\begin{equation}
\mathcal M(s,t)= \sum_{q=1}^N f_p(t-\Delta\tau(s,\vrho_q(s))),
\label{eq:Mdef}
\end{equation}
and $\Delta \tau$ given in (\ref{eq:DeltaTau}). Thus, the traces are
approximately, up to a multiplicative constant, equal to the matrix
$\{\cM(s_j,t_\ell)\}$, for $j = -n/2, \ldots, N/2$ and $\ell = -m/2, \ldots,
m/2$. We neglect the multiplicative constant hereafter, and call $\cM$
the matrix of traces.

We could take any pulse $f_p$ in (\ref{eq:Mdef}), but to simplify the
calculations in the remainder of the section, we assume it Gaussian,
modulated by a cosine at the central frequency,
\begin{equation}
f_p(t)=\cos(\omega_o t)e^{-B^2t^2/2}.
\label{eq:compPulse}
\end{equation}

\subsection{The annihilation filter}
\label{sect:filtanalysis}
Recall that the filter $\mathcal Q^{\vrho^e}$ defined by
(\ref{eq:OP1}) is a linear operator. Because we model the traces
(\ref{eq:Mdef}) as a superposition of pulses $f_p$ delayed by the
travel times to the targets, we can study the effect of $\mathcal
Q^{\vrho^e}$ on one trace at a time. 

Consider the trace from a stationary or moving target at location
$\vrho_j(s)$, for some integer $j$ satisfying $1 \le j \le N$, and
denote it by
\begin{equation}
\cM^j(s,t) = f_p(t-\Delta\tau(s,\vrho_j(s))). 
\end{equation}
The transformation $\mathbb{T}_+^{\vrho^e}$ defined in (\ref{eq:TT+})
maps the trace to
\begin{equation}
\left[\mathbb{T}_+^{\vrho^e} \cM^j\right](s,t) = 
f_p\left( t + \Delta \tau(s,\vrho^e)-\Delta \tau(s,\vrho_js))\right),
\label{eq:T1}
\end{equation}
and then we filter it with the difference operator $\mathbb{D}_s$, 
\begin{equation}
\mathbb{D}_s \left[\mathbb{T}_+^{\vrho^e} \cM^j\right](s,t) \approx
A(s,\vrho_j(s),\vrho^e) f'_p\left( t + \Delta \tau(s,\vrho^e)-\Delta
\tau(s,\vrho_j(s))\right).
\label{eq:T2}
\end{equation}
Here we approximate the finite difference by the derivative in $s$,
denote by $f_p'$ the derivative of the pulse, and let $A$ be the
``annihilation factor''
\begin{equation}
A(s,\vrho_j(s),\vrho^e)=\frac{d}{d
  s}\left[\tau(s,\vrho_j(s))-\tau(s,\vrho^e)\right].
\label{eq:T3}
\end{equation}

Using the scaling assumptions (\ref{eq:SC1}-\ref{eq:SC2}), we obtain
after a straightforward calculation that 
\begin{eqnarray}
\frac{d}{d s}\left[\tau(s,\vrho_j(s))-\tau(s,\vrho^e)\right] =
\frac{2}{c}\left[-\vu\cdot\vm_e+\frac{(V\vt-\vu)\cdot\Pp_e\Delta\vrho_j}{L}-
  s\frac{(2V\vt-\vu)\cdot\Pp_e\vu}{L}\right] + \nonumber \\ O\left(
\frac{a u_\perp}{c L}\right) + O\left( \frac{a VR^\cI}{c L^2}\right), \quad 
\label{eq:Dtau}
\end{eqnarray}
where 
\[
\Delta \vrho_j = \vrho_j(0) - \vrho^e, \quad 
\vm_e = \frac{\vr(0)-\vrho^e}{|\vr(0)-\vrho^e|}, \quad 
\mathbb{P}_e = I - \vm_e \vm_e^T.
\]
The annihilation factor of a stationary target is given by
\begin{equation}
\label{eq:AnnihSt}
A(s,\vrho_j,\vrho^e) = \frac{2V}{c} \frac{\vt \cdot \mathbb{P}_e \Delta 
\vrho_j}{L} +  O\left( \frac{a VR^\cI}{c L^2}\right).
\end{equation}
It is approximately linear in the cross-range component of the error
$\Delta \vrho^e$ of the estimated location of the target. For a moving
target the factor is
\begin{equation}
\label{eq:AnnihMv}
A(s,\vrho_j(s),\vrho^e) \approx \frac{2V}{c} \left[\frac{\vt \cdot
    \mathbb{P}_e \Delta \vrho_j}{L} - \frac{u}{V} - 2 \frac{Vs}{L}
  \frac{u_\perp}{V}\right] + O\left( \frac{a u_\perp}{c L}\right) +
O\left( \frac{a VR^\cI}{c L^2}\right).
\end{equation}
It is much larger than (\ref{eq:AnnihSt}) for stationary targets that
are in the vicinity of $\vrho^e$, and satisfy the condition
\begin{equation}
|\Delta \vrho_j| \ll L \left[ O \left(\frac{u}{V}\right) + O\left( 
\frac{a}{L} \frac{u_\perp}{V} \right) \right].
\label{eq:AnnihCond}
\end{equation}
Consequently, the filter $\mathcal Q^{\vrho^e}$ defined by
(\ref{eq:OP1}) emphasizes the traces from the moving targets and
suppresses the traces from the stationary targets near $\vrho^e$, as
stated in section \ref{sect:ANNIH}.

\begin{remark}
In the GOTCHA regime, and for one degree apertures as described in
section \ref{sect:scales}, the cross-range resolution of a preliminary
image is approximately $2.5$m. Thus, we can assume that for the target
that we wish to annihilate $|\Delta \vrho_j| \approx 2.5$m, and
condition (\ref{eq:AnnihCond}) is satisfied for targets moving at
range speed that exceeds $0.2$m/s.
\end{remark}

Note also from (\ref{eq:Dtau}), with $\vrho^e$ replaced by $\vrho_o$,
that the slope of the trace is approximately given by the range speed
$u = \vu \cdot \vm_o$, as assumed in section \ref{sect:rotationopt},
as long as
\[
u \gg V \frac{\vt \cdot \mathbb{P}_o \Delta \vrho_j}{L} \sim 
V \frac{R^{\cI}}{L}.
\]
If the target is a slower mover, than its trace looks similar to that
of a stationary target and it cannot be separated. The curvature of
the trace is determined by the cross-range velocity of the target, as
seen from the third term in the right hand side of equation
(\ref{eq:Dtau}).

% ------------------------------------------
\subsection{Analysis of data separation with robust PCA}
\label{sect:rpcaanalysis}
In this section we review briefly from \cite{BCP12} the analysis of
the rank of the matrix of traces for simple scenes. The goal is to
show that indeed, traces from stationary targets should be expected to
give a low rank contribution whereas those from the moving targets
should give a high rank contribution. This is the key assumption in
the data separation approach described in section \ref{sect:RPCA}.

\subsubsection{One target}
\label{sect:single}
We begin with the case of one target at location $\vrho(s) =
(\brho(s),0)$.  It moves at constant speed $\vu=(\bu,0)$ over the slow
time interval $|s|\le S(a)$, so we can write that
\begin{equation}
\vrho(s) = \vrho + s \, \vec{\bu}, \quad |s| \le S(a),
\label{eq:MOD9}
\end{equation}
where
\[
\vrho = \vrho(0).
\]
The matrix denoted by $\cM$ is defined by the discrete samples of
\begin{equation}
  \mathcal M(s,t) = \cos\left[\om_o(t - \Delta\tau(s,\vrho(s)))\right]
  \exp \left[ -\frac{B^2}{2} (t - \Delta\tau(s,\vrho(s)))^2 \right],
\label{eq:MOD8}
\end{equation}
where $\Delta \tau$ is defined by (\ref{eq:DeltaTau}), and we used the
pulse (\ref{eq:compPulse}).  The rank of $\cM$ is the same as that of
the symmetric, square matrix $\mathcal Y \in \mathbb{R}^{(n+1)\times
  (n+1)}$, given by
\begin{equation}
  \mathcal Y = \cM \cM^T. 
\label{eq:covar1}
\end{equation}
Its entries are the discrete samples of the function 
\begin{equation}
  \mathcal Y(s,s') = \sum_{q=-m/2}^{m/2} \cM(s,t_q) \cM(s',t_q) \approx
  \frac{1}{\Delta t}\int_{-\infty}^\infty dt\ \cM(s,t)\cM(s',t).
\label{eq:covar1p}
\end{equation}
We assume here that $\Delta t$ is small enough to approximate the
Riemann sum over $q$ by the integral over $t$. Since the traces vanish
for $|t| > \Delta s/2$, we can extend the integral to the whole real
line.

We showed in \cite{BCP12} that under our scaling assumptions, 
the function (\ref{eq:covar1p}) takes the form
\begin{equation}
  \mathcal Y(s,s') \approx \frac{\sqrt{\pi}}{2 B \Delta t} \cos \left[ \om_o
    \alpha (s-s')\right] \exp \left[-\frac{(B \alpha)^2(s-s')^2}{4}
    \right],
\label{eq:Toep1}
\end{equation} 
with dimensionless parameter
\begin{equation}
  \alpha =\frac{2 \vu \cdot \vm_o}{c}-\frac{2V \vt \cdot \Pp_o(
\vrho-\vrho_o)}{
    cL}+\frac{2 \vu \cdot \Pp_o(\vrho-\vrho_o)}{cL}.
\label{eq:Toep2}
\end{equation}
That is to say, the matrix (\ref{eq:covar1}) is approximately
Toeplitz, with diagonals given by the entries in the sequence
$\{y_j\}_{j \in \mathbb{Z}}$, 
\[
\mathcal Y_{j\ell} \approx y_{j-\ell},
\]
where 
\begin{equation}
  y_j =\frac{\sqrt{\pi}}{2B\Delta t}e^{-\frac{(B \Delta s |\alpha| 
      j)^2}{4}}\cos(\om_o \Delta s |\alpha|   \, j).
\label{eq:defcj}
\end{equation}
There are many slow time samples in an aperture $(n\gg 1)$, so the
matrix $\mathcal Y$ is also large, and we can approximate its rank using the
asymptotic Szeg\"{o} theory \cite{szego,bottcher}. 

Note that multiplication of a large Toeplitz matrix with a vector can
be interpreted approximately as a convolution, which is why the
spectrum can be estimated from the symbol $\hat y(\theta)$, the
Fourier transform of (\ref{eq:defcj}),
\begin{equation}
  \hat y(\theta)= \sum_{j=-\infty}^{\infty} y_j e^{ij\theta},
  \quad\quad \theta\in(-\pi,\pi).
\label{eq:series}
\end{equation}
Szeg\H{o}'s first limit theorem \cite{bottcher} gives that 
\begin{equation}
\label{eq:RANK_AS}
  \lim_{n\to\infty} \frac{\cN(n;\beta_1,\beta_2)}{n+1}=
  \frac{1}{2\pi}\int_{-\pi}^{\pi}1_{[\beta_1,\beta_2]}\left(\hat
  y(\theta)\right) d\theta,
\end{equation}
where $1_{[\beta_1,\beta_2]}$ is the indicator function of the
interval $[\beta_1,\beta_2]$ and $\cN(n;\beta_1,\beta_2)$ is the
number of eigenvalues of $\mathcal Y$ that lie in this interval. Using this 
result, we estimated the rank in \cite{BCP12} as  
\begin{equation}
  \mbox{rank}\left[\mathcal Y\right] := \cN\left(n;\epsilon \|\hat
  y\|_{\infty},\infty\right), \quad\quad \|\hat y\|_{\infty} =
  \sup_{\theta \in (-\pi,\pi)} |\hat y(\theta)|.
\end{equation}
Here $0< \epsilon \ll 1$ is a small threshold parameter, and $\|\hat
y\|_{\infty}$ is of the order of the largest singular value of
$\mathcal Y$.  The Szeg\H{o} theory \cite{szego,bottcher} gives that
this singular value is approximated by the maximum of the symbol.
The right hand side in (\ref{eq:RANK_AS}) is calculated explicitly in 
\cite{BCP12}, and the rank estimate is 
\begin{eqnarray}
  \frac{\mbox{rank}\left[\mathcal Y\right]}{n+1}&\approx&
  \frac{1}{2\pi}\int_{-\pi}^{\pi} 1_{[ \epsilon \|\hat y\|_{\infty},
      \infty)} \left(\hat y(\theta)\right) d\theta \nonumber \\ &=&
    \min\left(\frac{2|\alpha| B\Delta s\sqrt{\log
        1/\epsilon}}{\pi},1\right).
\label{eq:AS_RANK}
\end{eqnarray}
It grows linearly with $\alpha$, meaning that the larger the offset
$\vrho-\vrho_o$ in the cross-range direction is, and the faster the
target moves, the higher the rank.

\begin{figure}[t]
\centering
\includegraphics[width=.45\columnwidth]{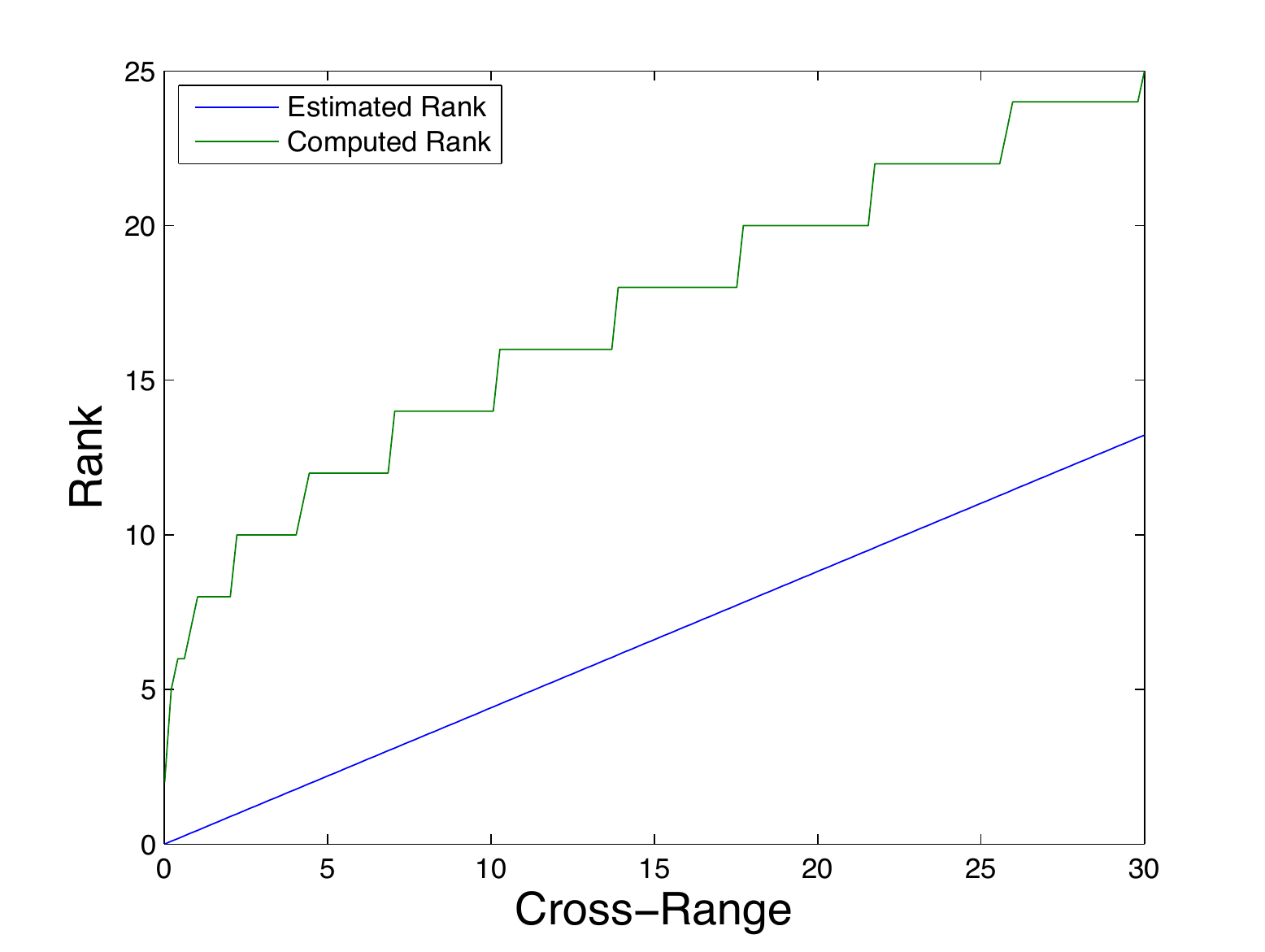}
\includegraphics[width=.45\columnwidth]{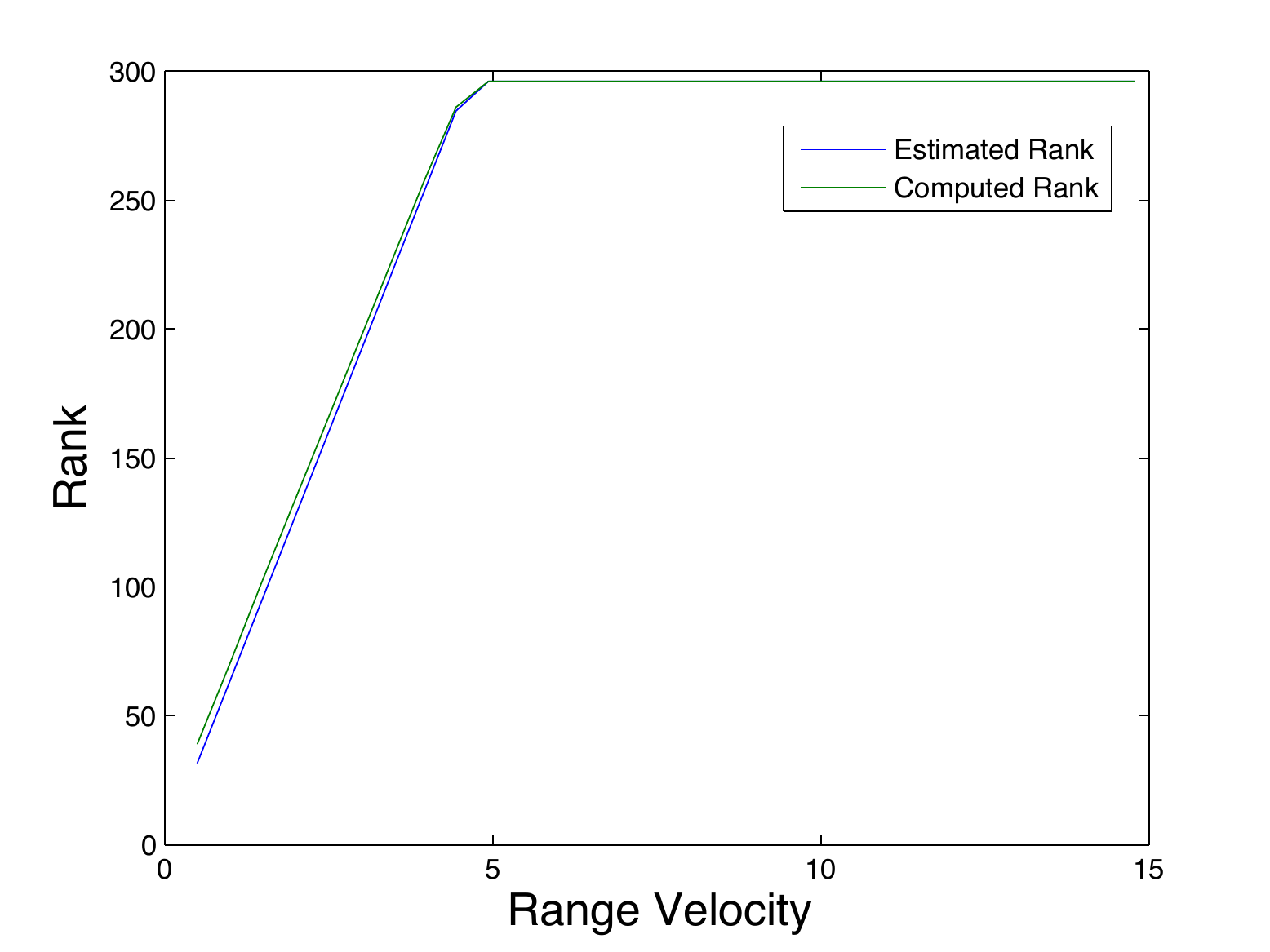}
\caption{Comparison of the computed and estimated rank of $\mathcal Y$ for a
single target. Left: Stationary target at various cross-range positions.
Right: Moving target with various velocities.}
\label{fig:rank1}
\end{figure}

\begin{figure}
\centering
\includegraphics[width=.5\columnwidth]{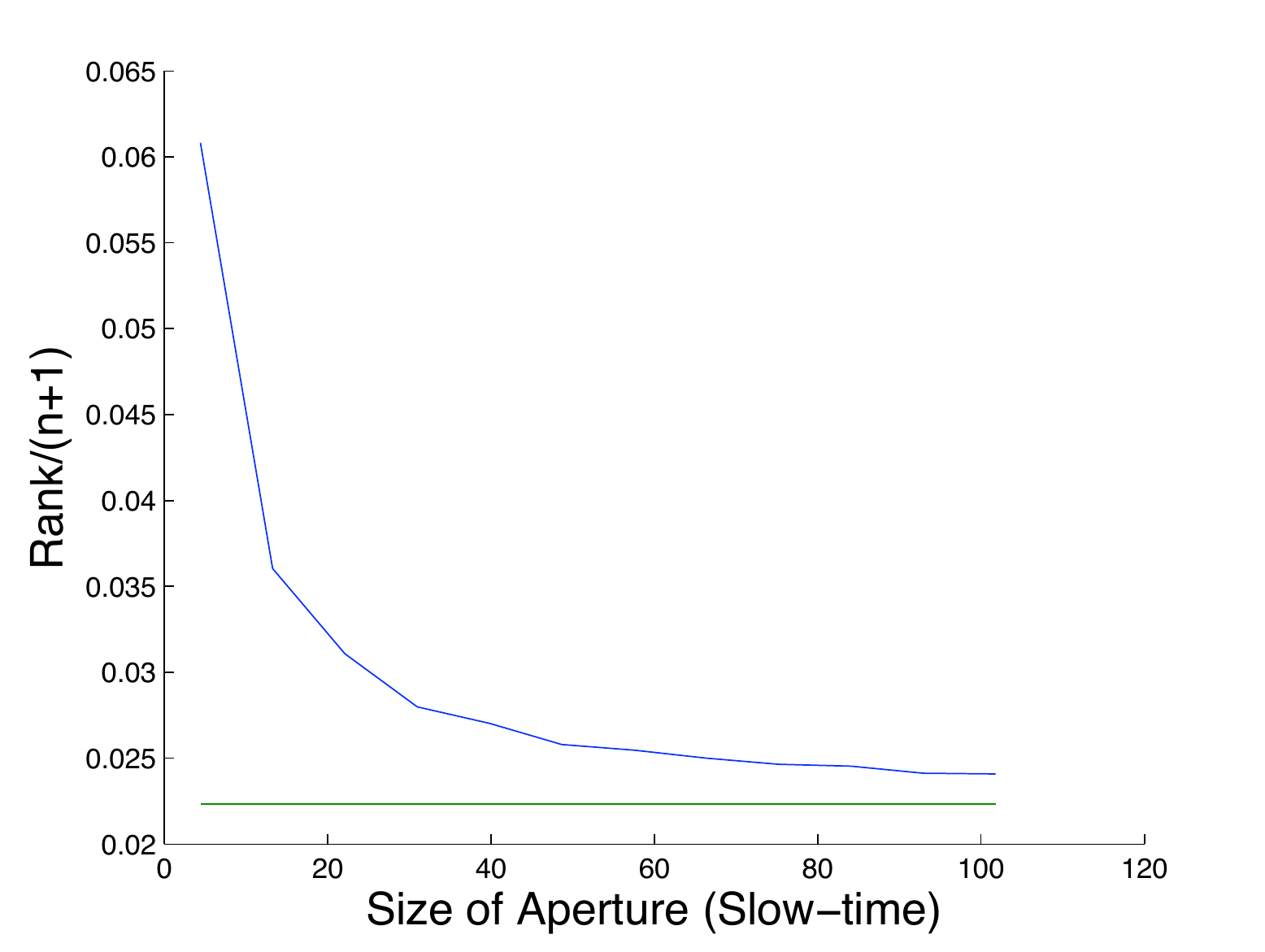}
\caption{Convergence of the rank of $\mathcal Y$ normalized by the
  size $(n+1)$. The blue line is the computed value and the green line
  is the asymptotic estimate. The slow time sampling is fixed, so the
  longer the aperture is, the larger $n$.}
\label{fig:convergeRank}
\end{figure}

We show in Figure \ref{fig:rank1} an illustration of this result.  The
left plot is for a stationary target and the right plot is for a
moving target. The green line shows the numerically computed rank of
the matrix $\mathcal Y$, while the blue line shows its asymptotic
estimate (\ref{eq:AS_RANK}). Note that both curves show the same
growth.  In the case of the stationary target there is a mismatch of
the magnitude of the computed and estimated rank. This is because $n$
is not large enough in the simulation.  A good approximation is
obtained when $y_j \approx 0$ for $j \approx n$, because then we can
approximate the series (\ref{eq:series}) that defines the symbol by
the truncated sum for indices $|j| \le n$. This is not the case in
this simulation, so there is a discrepancy in the estimated
rank. However, if we increase the aperture, thus increasing $n$, the
result improves. Figure \ref{fig:convergeRank} illustrates this fact
by showing how the rank normalized by the size of the matrix has the
predicted asymptote as $n$ increases.

In the right plot of Figure \ref{fig:rank1} we compare the computed
rank (in green) and the asymptotic estimate (\ref{eq:AS_RANK}) (in
blue) for a moving target located at $\vrho = \vrho_o$ at $ s= 0$. The
rank is plotted as a function of the range component of the velocity.
Note that, as expected, the rank increases with the velocity, and
therefore with $|\alpha|$. For the moving target, the asymptotic
estimate is in closer agreement to the computed one.  This is because
in this case the entries in the sequence $\{y_j\}_{j \in \mathbb{Z}}$
decay faster with $j$.  Finally, we can observe that the rank is much
larger for the moving target than the stationary target, even for
small velocities.

\subsubsection{Two targets}
\label{sect:double}
When there are two targets in the scene, we obtain from model
(\ref{eq:Mdef}) and the choice (\ref{eq:compPulse}) of the pulse that 
\begin{equation}
  \mathcal M(s,t) = \sum_{j=1}^2\cos\left[\om_o(t -
    \Delta\tau(s,\vrho_j))\right] \exp \left[ -\frac{B^2}{2} (t -
    \Delta\tau(s,\vrho_j))^2 \right].
\label{eq:TWO1}
\end{equation}
We restrict the discussion to the case of two stationary targets. The
extension to moving targets is straightforward. It amounts to
adjusting the parameters $\alpha_1$ and $\alpha_2$ defined below in
(\ref{eq:alphaj}), by adding two linear terms in the target velocity,
as in equation (\ref{eq:Toep2}).

We obtain after a calculation given in detail in \cite{BCP12} that 
the entries of the matrix (\ref{eq:covar1}) are the discrete slow time 
samples  of the function 
\begin{align}
\mathcal Y(s,s')&\approx\frac{\sqrt{\pi}}{2B\Delta t}\left\{ \sum_{j=1}^2
\cos[\omega_o\alpha_j(s-s')]\exp \left[-\frac{(B\alpha_j)^2
(s-s')^2}{4}\right] \right.\nonumber\\ &+ \cos[\omega_o(\alpha_1
s-\alpha_2 s'+\beta)] \exp\left[-\frac{B^2(\alpha_1 s-\alpha_2
s'+\beta)^2}{4}\right]\label{eq:covar2}\\ &\left. +
\cos[\omega_o(\alpha_1 s'-\alpha_2 s+\beta )]
\exp\left[-\frac{B^2(\alpha_1 s'-\alpha_{2}
s+\beta)^2}{4}\right]\right\},\nonumber
\end{align}
where 
\begin{equation}
\alpha_j =-\frac{2V \vt \cdot \mathbb{P}_o
(\vrho_j-\vrho_o)}{cL}, \quad j = 1,2 ,
\label{eq:alphaj}
\end{equation}
and 
\begin{equation}
\beta = \frac{2}{c}\sum_{j=1}^2 (-1)^{j} \left\{ \vm_o \cdot
(\vrho_j-\vrho_o) + \frac{[\vm_o \cdot (\vrho_j-\vrho_o)]^2}{2 L}
\right\}.
\label{eq:beta}
\end{equation}
Thus, the matrix $\mathcal Y$ has a more complicated structure. The
first term in (\ref{eq:covar2}), the sum over $j$, gives rise to a
Toeplitz matrix, as before. The last two terms give matrices that are
approximately g-Toeplitz or g-Hankel, depending on the sign of the
ratio $\alpha_2/\alpha_1$.

\begin{definition}
An $(n+1) \times (n+1)$ g-Hankel matrix  $\mathcal H$ with shift $g \in
\mathbb{Z}^{+}$ is defined by a sequence $\{h_j\}_{j \in \mathbb{N}}$ as 
\[
\mathcal H_{j\ell} = h_{j-1+g(\ell-1)}, \qquad j, \ell = 1,\ldots, n+1.
\]
The matrix is Hankel when $g = 1$.  A g-Toeplitz matrix is defined
similarly, by replacing $g$ with $-g$.
\end{definition}

It is easier to analyze the spectrum of the matrix $\mathcal Y$, when 
the reference point $\vrho_o$ is chosen so that 
\begin{equation}
\frac{\alpha_2}{\alpha_1} < 0, \quad \mbox{and}
\quad g:=\left|\frac{\alpha_2}{\alpha_1}\right| \in \mathbb{N}.
\label{eq:assumerho_o}
\end{equation}
Then, $\mathcal Y$ is given by 
\begin{equation}
\mathcal Y = \mathcal T + \mathcal H + \mathcal H^T,\label{eq:sumTH}
\end{equation}
with the Toeplitz and g-Hankel matrices 
\begin{equation}
\label{eq:Tg}
\mathcal T_{j\ell} = y_{j-\ell}, \qquad \mathcal H_{j\ell}=h_{(j-1)+g(\ell-1)},
\qquad j, \ell = 1, \ldots, n+1.
\end{equation}
The sequences $\{y_j\}_{j \in \mathbb{Z}}$ and $\{h_j\}_{j \in \mathbb{N}}$ 
that define these matrices are 
\begin{equation}
y_j = \frac{\sqrt{\pi}}{2B\Delta t} \sum_{\ell=1}^2 \cos \left( \om_o
\Delta s |\alpha_\ell| j \right) e^{-\frac{(B \Delta s |\alpha_\ell|j)^2}{4}},
\label{eq:seqc}
\end{equation}
and 
\begin{equation}
h_j = \frac{\sqrt{\pi}}{2B\Delta t}e^{-\frac{[B \Delta s |\alpha_1| (j
      + \zeta)]^2}{4}}\cos[\om_o \Delta s |\alpha_1| (j + \zeta)],
\qquad \zeta = \frac{\beta}{|\alpha_1| \Delta s}.
\label{eq:seqh}
\end{equation}

If the range offsets from the reference point $\vrho_o$ are
sufficiently large, so that
\[
B |\beta| = \frac{2 B}{c} \left|\sum_{\ell=1}^2 (-1)^\ell \{ \vm_o \cdot
(\vrho_\ell - \vrho_o) + \frac{[\vm_o \cdot (\vrho_\ell - \vrho_o)]^2}{2 L}
\right| \gg 1,
\]
then the g-Hankel matrix has small entries, and $\mathcal Y$ is
approximately Toeplitz, as in the single target case. The meaning of
this is that the difference of the travel times between the SAR
platform and such targets is larger than the compressed pulse width,
which makes their interaction in (\ref{eq:covar2}) negligible.
However, if the range offsets are small, the interaction plays a role
and the structure of the matrix $\mathcal Y$ is given by equation
(\ref{eq:sumTH}).  In this case, the estimate of the rank of $\mathcal
Y$ follows from the recent results in \cite{sesana,fasino,tilli}. They
state that the g-Hankel terms $\mathcal H + \mathcal H^T$ have a
negligible effect on the rank in the limit $n \to \infty$.  See
\cite{BCP12} for more details.  Thus, we can conclude in either case
that the rank estimate of $\mathcal Y $ is given by equation
(\ref{eq:RANK_AS}), in terms of the symbol $\hat y(\theta)$ defined by
(\ref{eq:series}), using the sequence $\{y_j\}_{j\in \mathbb{Z}}$ with
entries (\ref{eq:seqc}).  

\begin{figure}[t]
\centering
\includegraphics[width=.9\columnwidth]{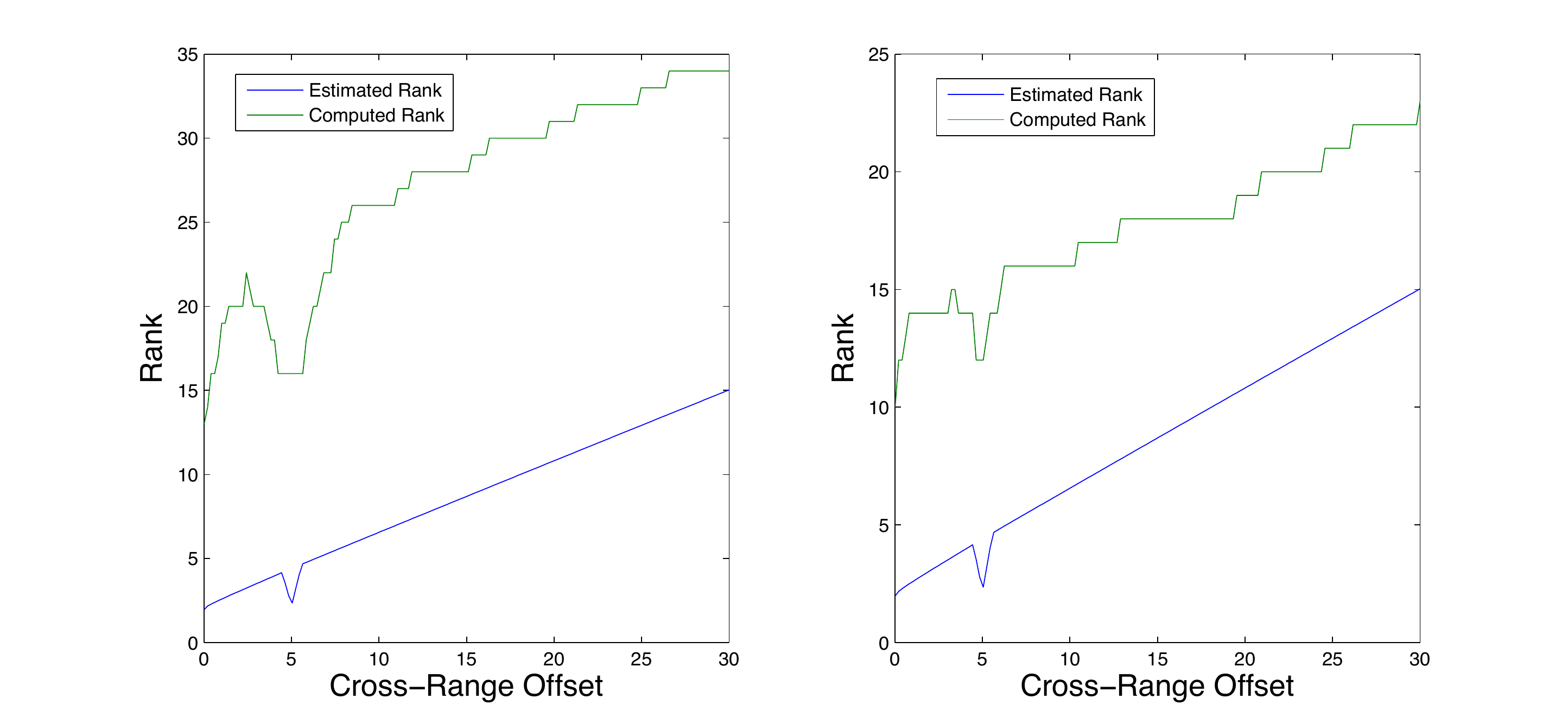}
\caption{Computed and estimated rank of matrix $\mathcal Y$ for a
  scene with two stationary targets. On the left, one target is fixed
  at location $\vrho_1 =(5,5,0)$m and the location of the other varies
  on the line segment between $(-5,0.01,0)$m and $(-5,30,0)$m.  On the
  right, one target is fixed at location $\vrho_1 =(0.15,5,0)$m and
  the location of the other varies on the line segment between
  $(-0.15,0.01,0)$m and $(-0.15,30,0)$m.}
\label{fig:eigenvals2}
\end{figure}
We plot in Figure \ref{fig:eigenvals2} the computed and estimated rank
for two different configurations of the stationary targets.  The
results on the left consider one target at location
$\vrho_1=(5,5,0)$m.  The second target location is varied between
$(-5,0.01,0)$m and $(-5,30,0)$m.  The range separation is large, equal
to $10$m, so the entries of the g-Hankel matrix $H$ are negligible
because $B|\beta|=41.47 \gg 1$.  The plot on the right considers one
target at $\vrho_1=(0.15,5,0)$m and a second target at location the
$\vrho_2$ varying between $(-0.15,0.01,0)$m and $(-0.15,30,0)$ m.  The
range separation is now small, and the g-Hankel matrix $\mathcal H$ is
no longer negligible.  Nevertheless, as predicted by the asymptotic
theory, the rank of matrix $\mathcal Y$ behaves essentially the same
as before.

\subsubsection{Conclusions of the analysis of rank}
\label{sect:rankDisc}
Several conclusions can be made from the analysis summarized above.
First, the traces from a moving target form a high rank matrix.  For a
fixed number $n$ of slow time samples, the rank grows linearly with
the target speed, until the matrix becomes full rank.  The traces from
a stationary target form a low rank matrix if the scene is small
enough. The rank grows linearly with the cross-range offset of the
target from the reference point $\vrho_o$, but the growth is at a
slower rate than for a moving target.

We also analyzed the rank for scenes with two stationary targets and
found that it increases linearly with the cross-range offsets of the
targets from $\vrho_o$. The expectation, confirmed by the numerical
simulations, is that the rank increases with the number of targets.
This is why the data separation with robust PCA should be done in
successive small time windows, with each window containing the traces
from only a few stationary targets. These traces give a matrix that is
low rank, and thus can be separated from the traces due to moving
targets. The results shown in Figure \ref{fig:rpcaEx2b} illustrate
this point.

% ------------------------------------------
% ------------------------------------------
\section{Numerical Simulations}
\label{sect:numeric}
We present here more numerical simulations that illustrate the
performance of the annihilation filters defined in sections
\ref{sect:ANNIH} and \ref{sect:rotationopt}. We have presented already
in sections \ref{sect:RPCA} and \ref{sect:sepdisc} three different
simulations for data separation with robust PCA. More results can be
found in \cite{BCP12}.

\subsection{Setup of the simulations}

We use the GOTCHA Volumetric SAR setup as described in section
\ref{sect:scales}.  The data traces are generated with the model
(\ref{eq:MOD1}).  In all the simulations, the point targets are
assumed to have identical reflectivity $\sigma_q=1$.  The motion
estimation results are obtained with the phase space algorithm
introduced in \cite{sar}.  This algorithm requires that we know the
location of the target at one moment during the measurement window.
We choose it at the center of the aperture, which is why there is no
error in the target trajectory at $s=0$.

For the robust PCA results, the principal component pursuit
optimization is solved with an augmented Lagrangian approach.  It
requires the computation of the top few singular values and
corresponding singular vectors of large and sparse matrices, which we
do with the software package PROPACK.

\subsection{Annihilation Filters}
\begin{figure}[!t]
\centering\includegraphics[width=.7\columnwidth]{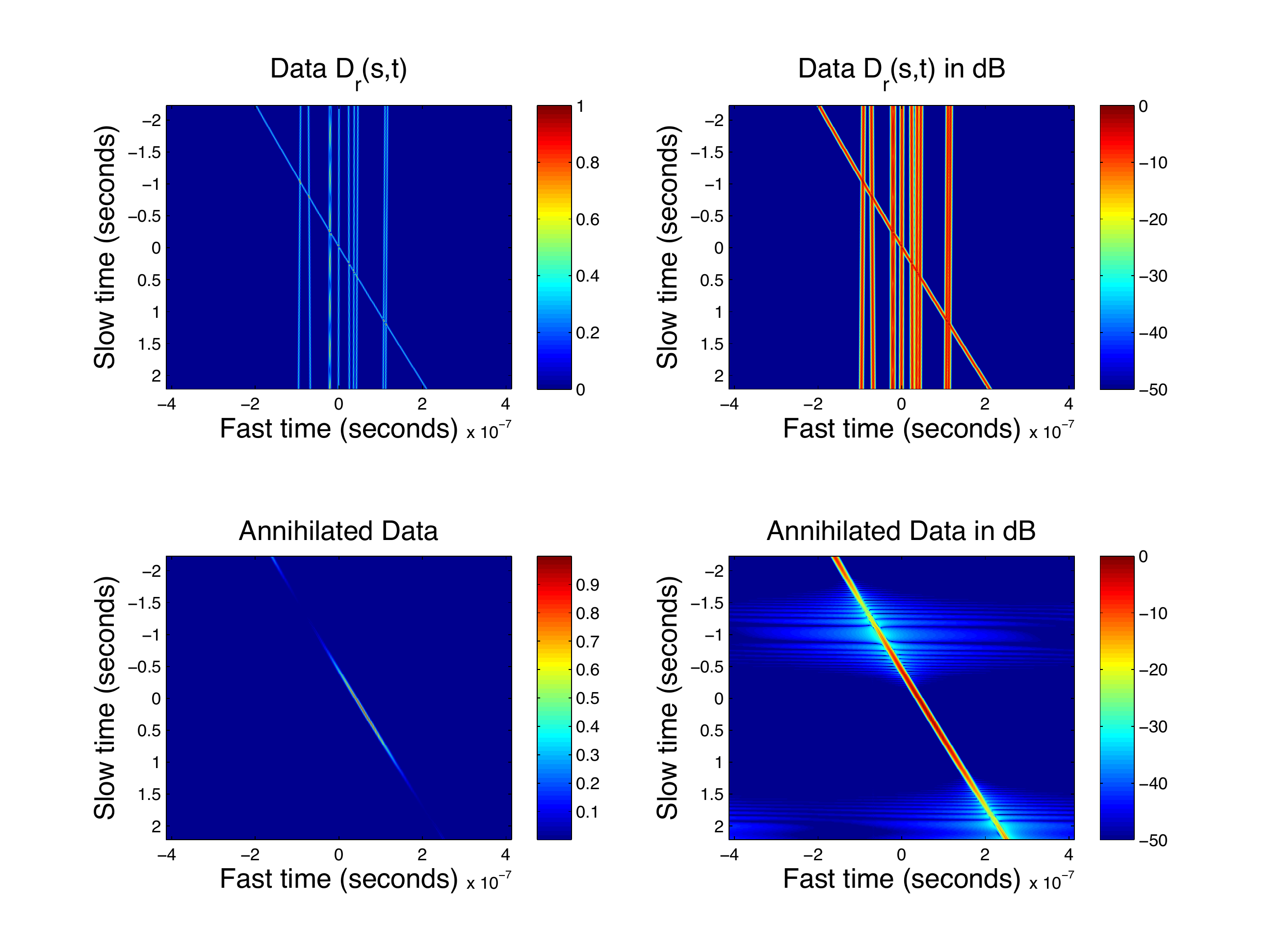}
\caption{Data separation with the annihilation filter for a scene with
  ten stationary targets and a moving one. The top row shows the data
  traces and the bottom row the trace from the moving target, as
  separated by the algorithm.  The right column is the left column
  plotted in dB scale.}
\label{fig:filter10scatt1ord}
\end{figure}
\begin{figure}[!t]
\centering\includegraphics[width=.7\columnwidth]{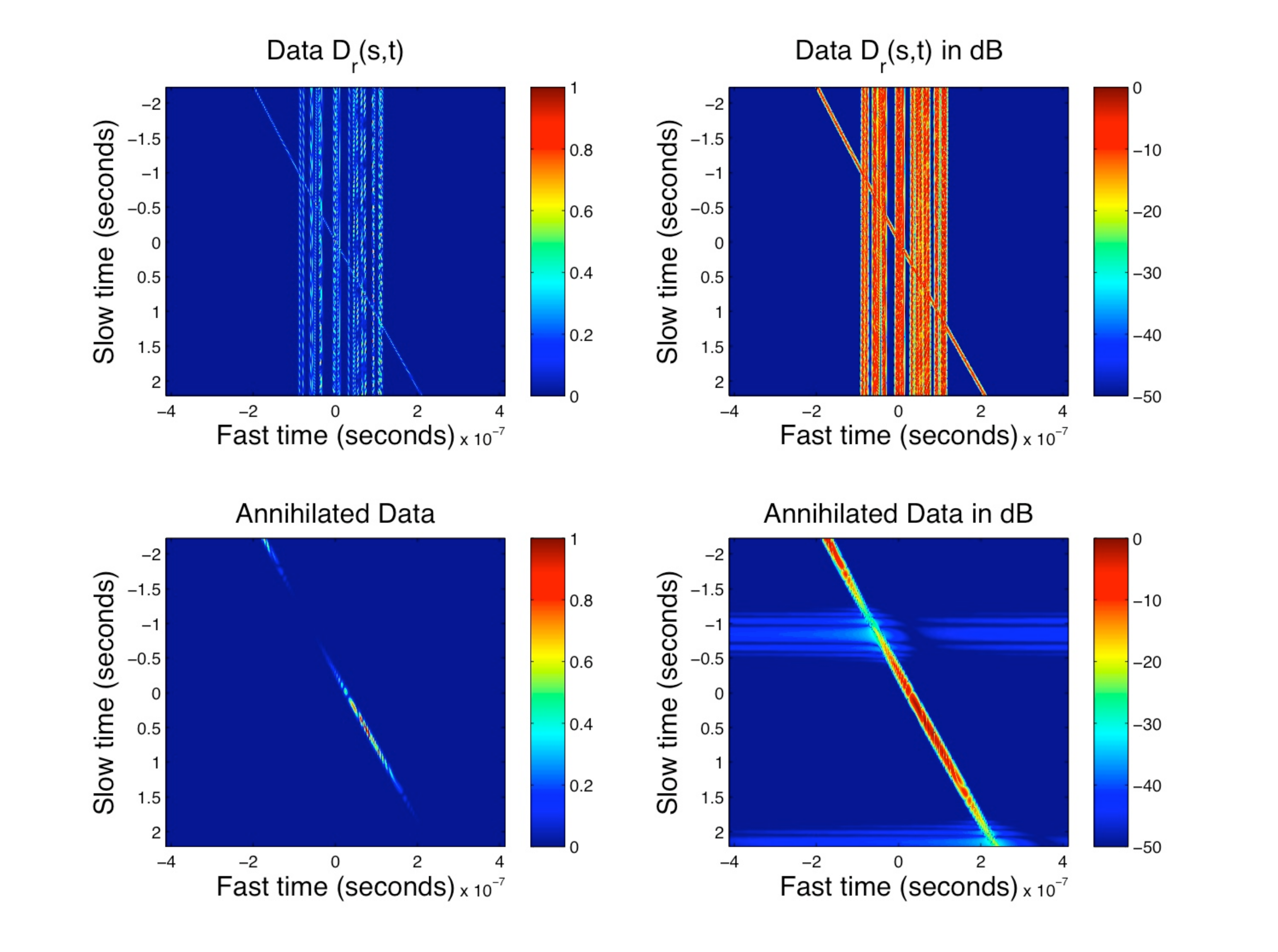}
\caption{Data separation with the annihilation filter for a scene with
  thirty stationary targets and a moving one. The top row shows the
  data traces and the bottom row the trace from the moving target, as
  separated by the algorithm.  The right column is the left column
  plotted in dB scale.}
\label{fig:filter30scatt1ord}
\end{figure}
\begin{figure}[!t]
\centering
\subfigure{\includegraphics[width=0.3\columnwidth]{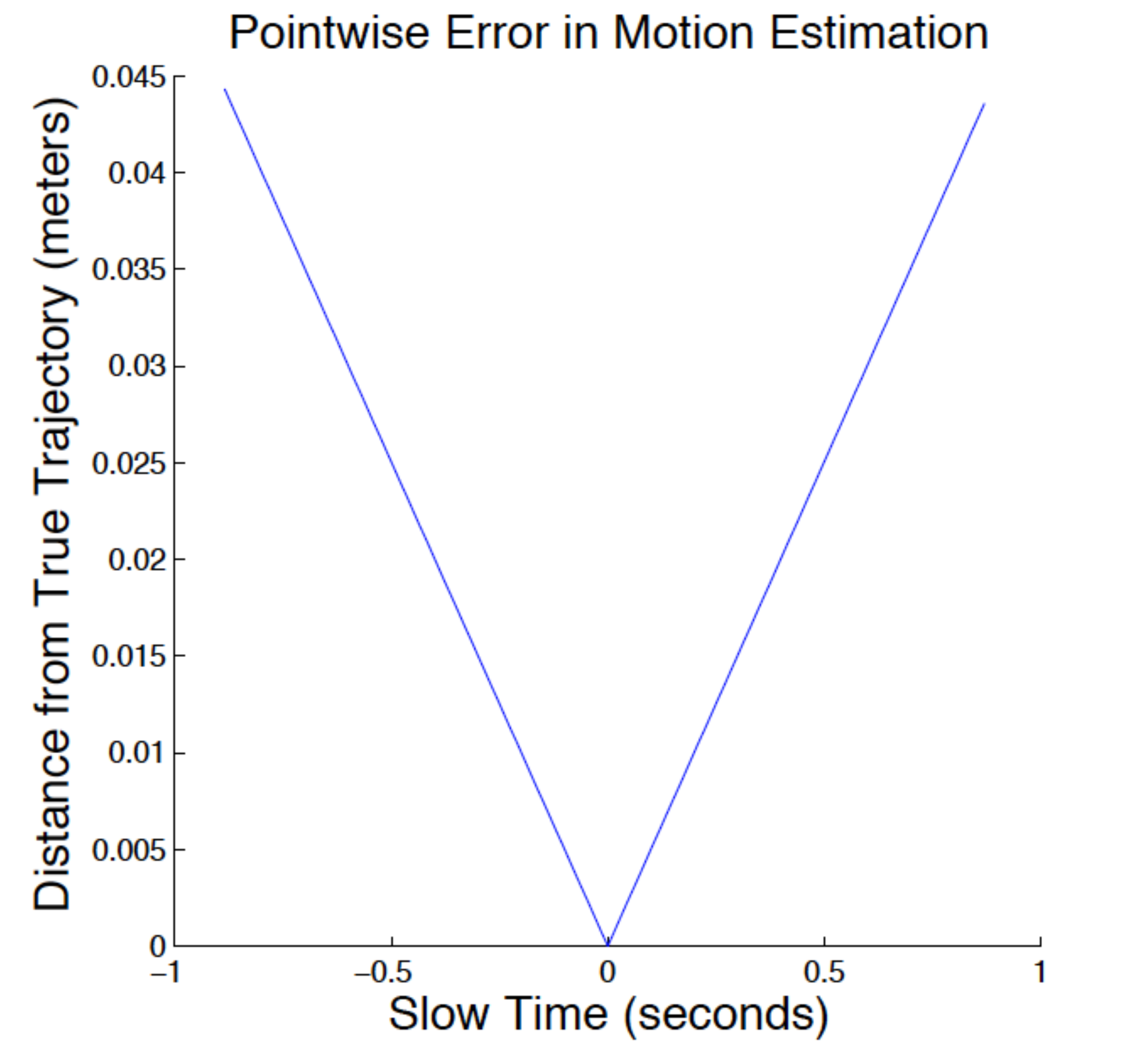}}
\subfigure{\includegraphics[width=0.285\columnwidth]{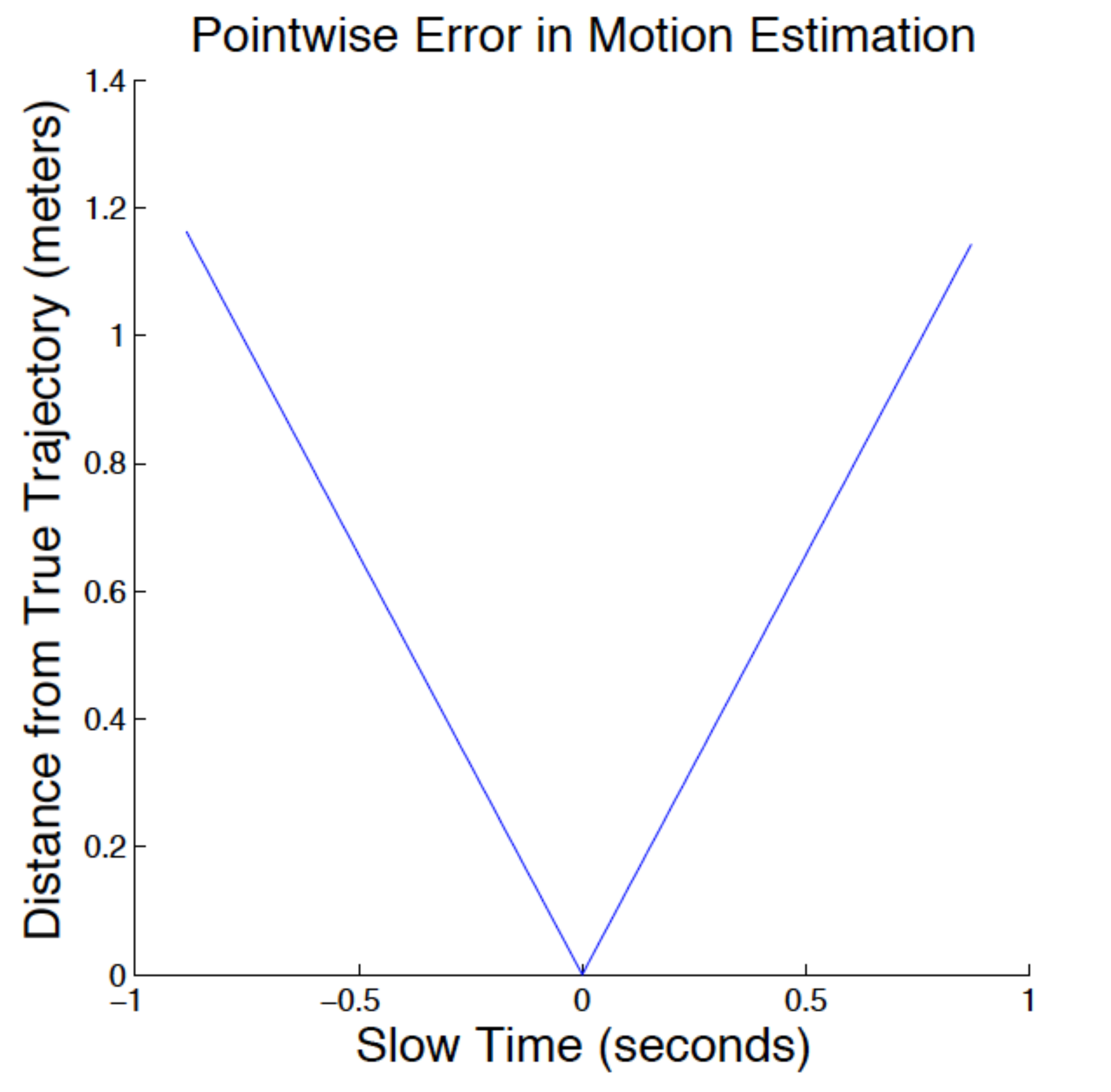}}
\caption{Error of the estimated target trajectory in a complex scene
  with ten stationary targets (left) and thirty stationary targets
  (right).  The motion estimation is done with the filtered data shown
  in the bottom row of Figures \ref{fig:filter10scatt1ord} and
  \ref{fig:filter30scatt1ord}.}
\label{fig:motEstscatt}
\end{figure}

We give two examples of data separation with the annihilation filters
described in section \ref{sect:ANNIH}.  In the first example, we have
a complex scene composed of ten stationary scatterers in an imaging
region of 50$\times$50$\mbox{m}^2$, and a moving target with velocity
$\vu=28/\sqrt{2}(1,1,0)$m/s.  We plot in the top row of Figure
\ref{fig:filter10scatt1ord} the data traces, and in the bottom row the
filtered traces.  The right column is the same as the left except it
is plotted in decibel (dB) scale to increase the contrast of the
figure.  The second example considers a scene with thirty stationary
targets and the same one moving target.  The results are in Figure
\ref{fig:filter30scatt1ord}.  In both examples we used the exact
stationary scatterer locations in the definition of the filters
(\ref{eq:OP1}). The errors in the estimated target
  trajectories, with the separated trace from the moving target for each example,
 are shown in Figure \ref{fig:motEstscatt}.

Note that the filtering and the motion estimation results are worse in
the second example, where there are more stationary targets to remove,
by applying repeatedly filters like (\ref{eq:OP1}) to the traces.
Computing many times numerical approximations of slow time derivatives
causes errors to accumulate. This is specially because we have high
frequency signals, and the sampling in slow time is limited by the
pulse repetition frequency of the SAR system.

\subsection{Data separation for Scene 1}
We already saw in section \ref{sect:sepdisc}, Figure
\ref{fig:rpca2mov}, the separation of the traces from the twenty
stationary and two moving targets in Scene 1.  Here we show the
further separation of the traces from the moving targets, as described
in section \ref{sect:rotationopt}.  The first moving target is located
at $(0,0,0)$m at $s = 0$, and is moving in the horizontal plane with
velocity ${\bf u} = \frac{28}{\sqrt{2}}(1,1)$m/s. The second target is
located at $(-5,5,0)$m at $s = 0$, and moves with velocity ${\bf u} =
\frac{14}{\sqrt{3}} \left(-1,\sqrt{2}\right)$m/s.

\begin{figure}[!t]
\centering
\includegraphics[width=0.9\columnwidth]{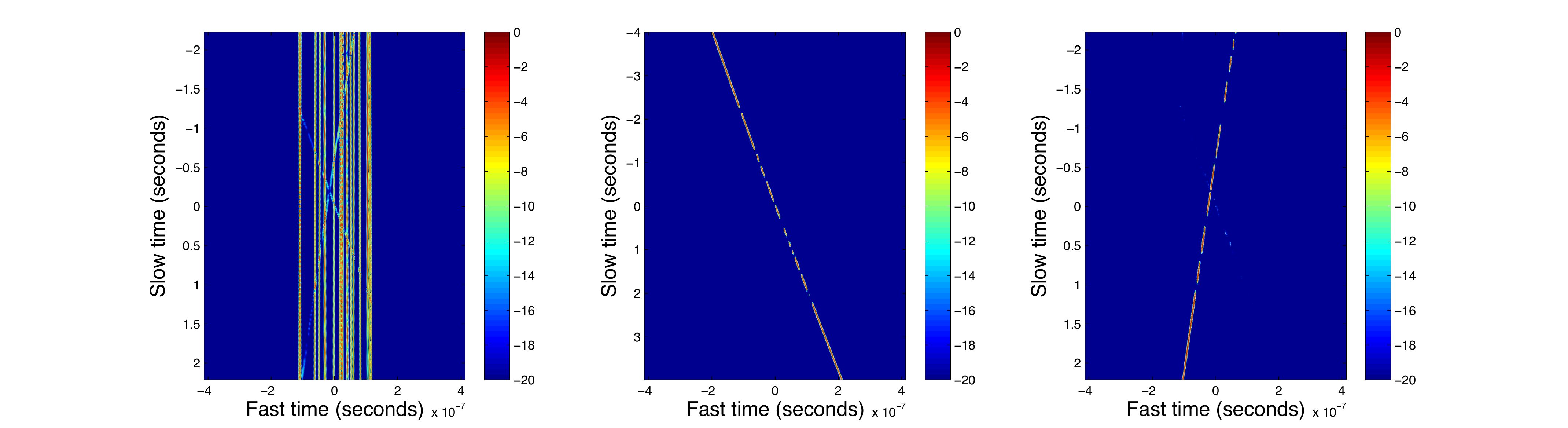}
\caption{Data separation for Scene 1. Left: traces from the stationary targets.
Middle and right: traces from the moving targets.}
\label{fig:rpca2movFinalSep}
\end{figure}

\begin{figure}[!t]
\centering
\subfigure{\includegraphics[width=0.33\columnwidth]{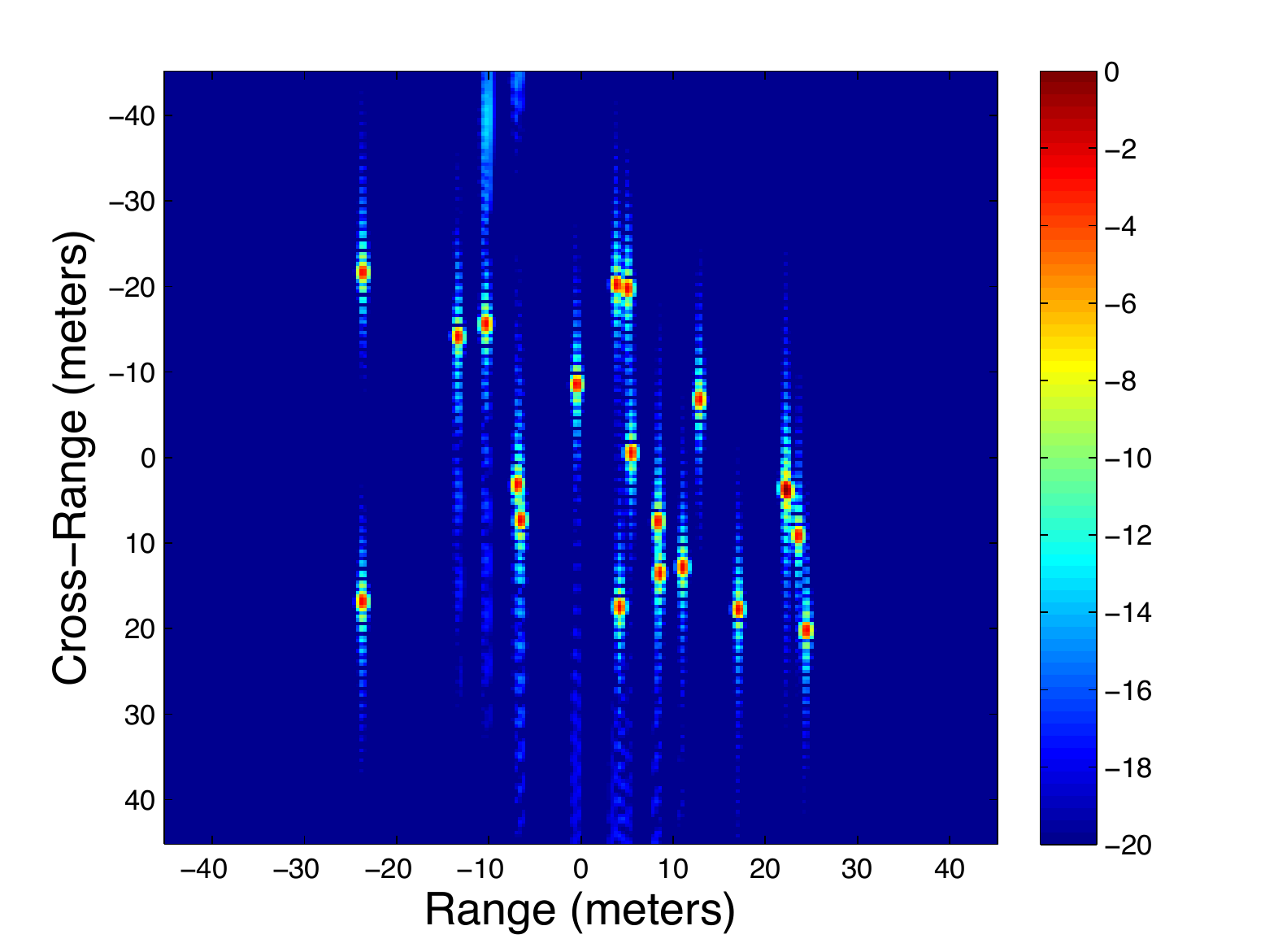}}\hspace{-.15in}
\subfigure{\includegraphics[width=0.33\columnwidth]{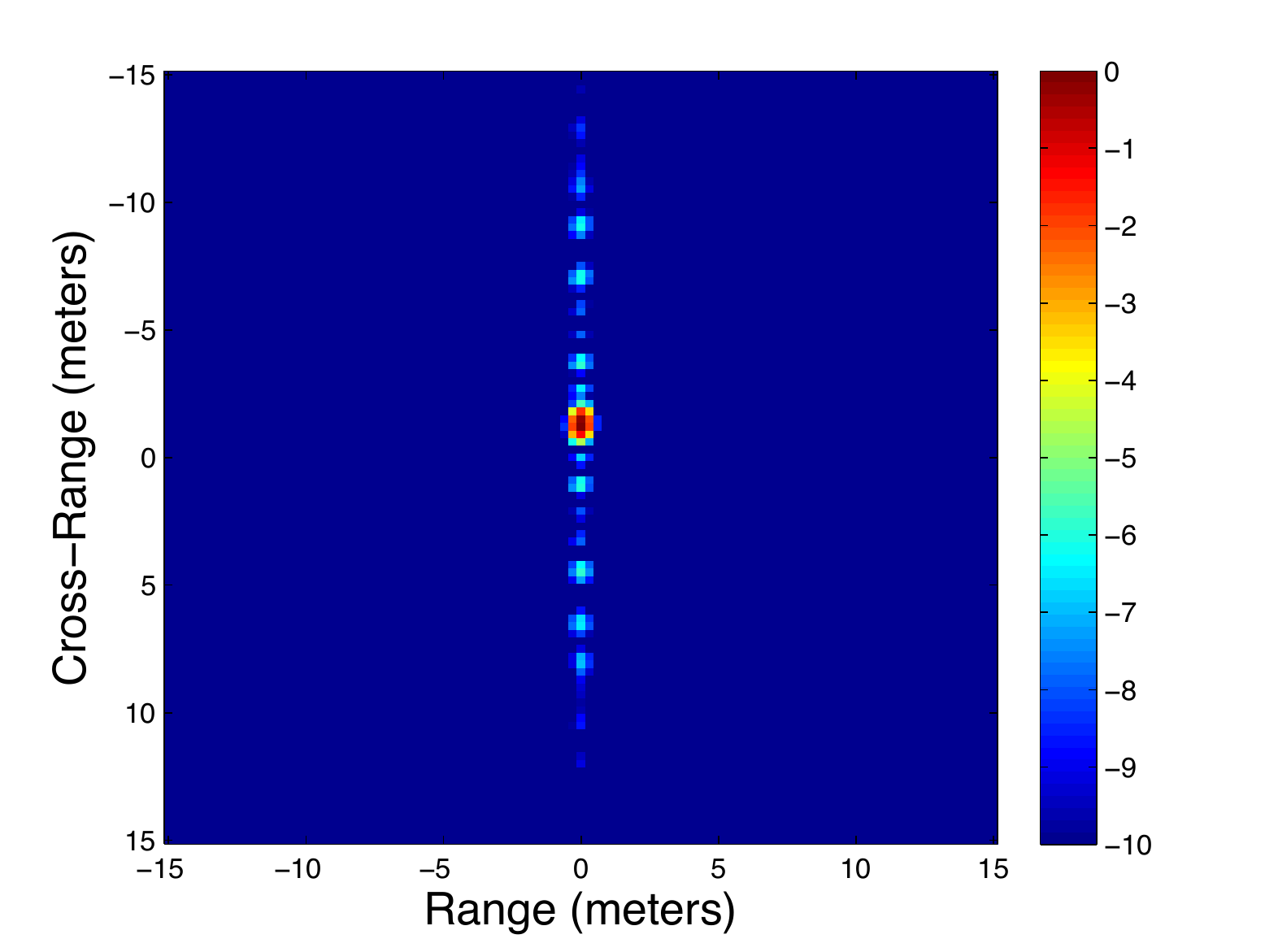}}\hspace{-.15in}
\subfigure{\includegraphics[width=0.33\columnwidth]{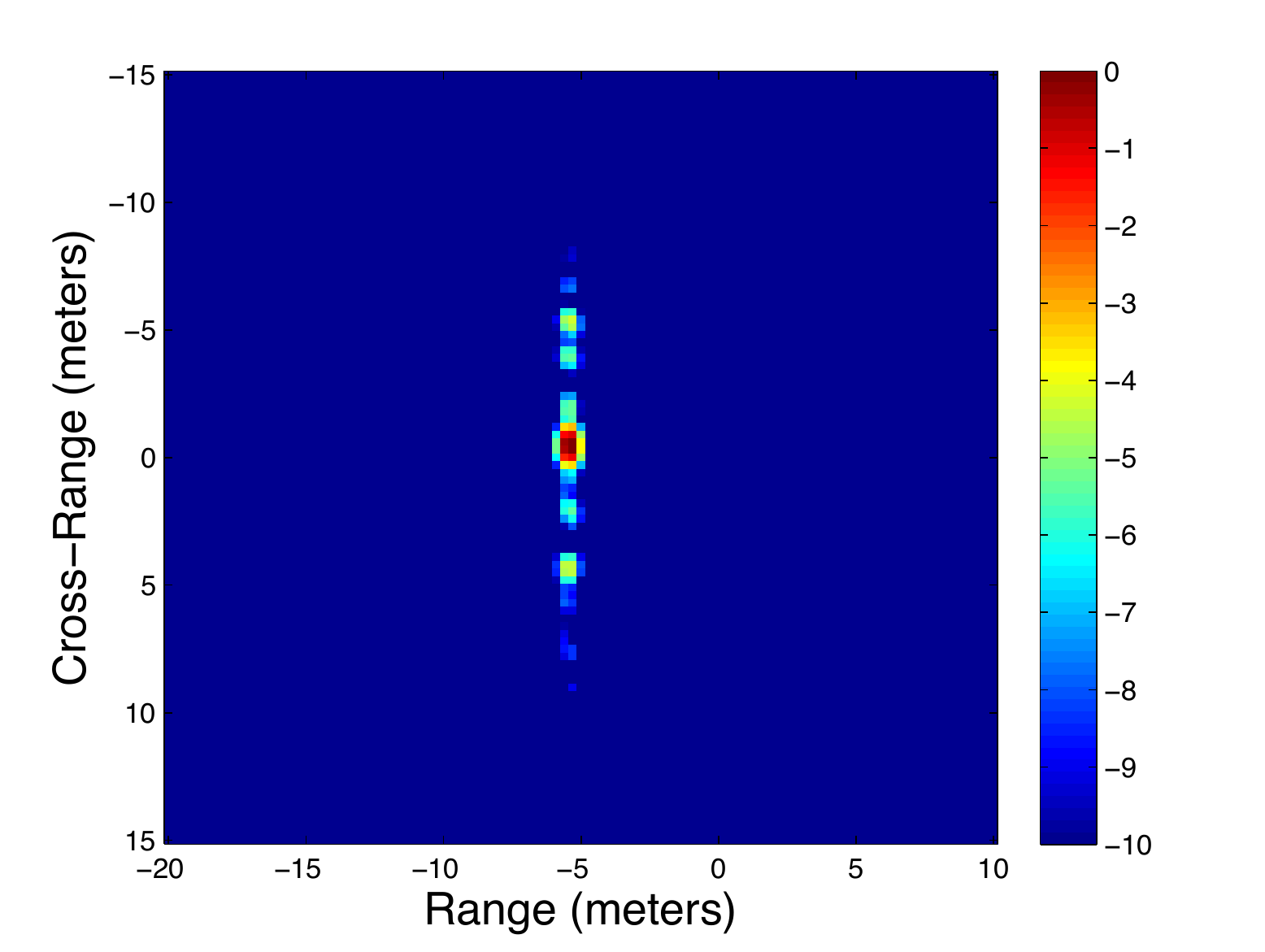}}
\caption{Images for Scene 1. Left: image with the traces from the
  stationary targets.  Middle and right: image with the traces the
  moving targets.}
\label{fig:sepImgs}
\end{figure}

We apply the method described in section \ref{sect:sepdisc} to the
traces shown in the right plot of Figure \ref{fig:rpca2mov}, separated
by robust PCA. We estimate the velocity of the first target,
apply the travel time transformation (\ref{eq:TTu+}) to the traces,
and then use the robust PCA approach to do the separation.  The output
of the data separation is shown in Figure \ref{fig:rpca2movFinalSep}.
The traces from the stationary targets, as obtained by the robust PCA
method are on the left. The traces from the two moving targets are in
the middle and right plots. The images obtained with these traces are
in Figure \ref{fig:sepImgs}. All the targets, stationary and moving
are now in focus.  The peak of the image for the second moving target 
is slightly displaced in cross-range due to the error in its velocity estimation
after separation.

% ------------------------------------------
\section{Summary}
\label{sect:conc}
We presented two approaches for separating the pulse and range
compressed echoes from stationary and moving targets in a complex
scene using a single synthetic aperture radar (SAR) antenna.  The
separation is required for estimating the motion of the targets and
focusing the images.

The first approach subtracts the echoes from the stationary targets
one by one, by performing the following four steps: (1) It estimates
the stationary target location $\vrho^e$ from a preliminary image. (2)
It transforms the data using a map $\mathbb{T}_+^{\vrho^e}$ so that
the echo from the target becomes approximately independent of the SAR
antenna location on the flight path.  (3) It annihilates the target
echo by exploiting this independence. (4) It undoes the transformation
at step (2) by applying the map $\mathbb{T}_{-}^{\vrho^e}$ to the
remaining echoes.

The second approach is based on the robust principle component
analysis (PCA) method, which is designed to decompose a matrix $\cM$
into its low rank part and high rank but sparse part. In our context,
the entries of the matrix $\cM$ are the samples of the pulse and range
compressed echoes measured by the SAR antenna.  The main observation
is that with appropriate pre-processing and windowing, the discrete
samples of the stationary target echoes form a low rank matrix, while
those from a few moving targets form a high rank sparse matrix. Thus,
they can be separated with robust PCA.

We presented a brief analysis of the two methods, and a numerical
comparison of their performance. Each method has its advantages and
disadvantages, but they complement each other. Thus, we can combine
them to improve the data separation results for complex and extended
imaging scenes.

% ------------------------------------------
\section*{Acknowledgement}

The work of L. Borcea was partially supported by the AFSOR Grant
FA9550-12-1-0117, by Air Force-SBIR FA8650-09-M-1523, the ONR Grant
N00014-12-1-0256, and by the NSF Grants DMS-0907746, DMS-0934594.  The
work of T. Callaghan was partially supported by Air Force-SBIR
FA8650-09-M-1523 and the NSF VIGRE grant DMS-0739420.
The work of G. Papanicolaou was supported in part by AFOSR grant FA9550-11-1-0266.

\bibliographystyle{plain} \bibliography{SBIR}

\begin{thebibliography}{10}

\bibitem{barbarossa1992detection}
S.~Barbarossa and A.~Farina.
\newblock {Detection and imaging of moving objects with synthetic aperture
  radar}.
\newblock {\em IEE Proceedings-F}, 139(1):79--88, 1992.

\bibitem{BCP12}
L.~Borcea, T.~Callaghan, and G.~Papanicolaou.
\newblock {Synthetic Aperture Radar Imaging and Motion Estimation via Robust
  Principal Component Analysis}.
\newblock {\em arXiv:1208.3700}, 2012.

\bibitem{sar}
L.~Borcea, T.~Callaghan, and G.~Papanicolaou.
\newblock {Synthetic Aperture Radar Imaging with Motion Estimation and
  Autofocus}.
\newblock {\em Inverse Problems}, 28:045006, 2012.

\bibitem{delcueto}
L.~Borcea, F.~Gonzalez~del Cueto, G.~Papanicolaou, and C.~Tsogka.
\newblock {Filtering deterministic layering effects in imaging}.
\newblock {\em SIAM Multiscale Model. Simul.}, 7(3):1267--1301, 2009.

\bibitem{delcueto2}
L.~Borcea, F.~Gonzalez~del Cueto, G.~Papanicolaou, and C.~Tsogka.
\newblock {Filtering random layering effects in imaging}.
\newblock {\em SIAM Multiscale Model. Simul.}, 8(3):751--781, 2010.

\bibitem{bottcher}
A.~B\"{o}ttcher and B.~Silbermann.
\newblock {\em Introduction to Large Truncated Toeplitz Matrices}.
\newblock Springer, 1999.

\bibitem{candesRPCA}
E.~J. Cand\`{e}s, X.~Li, Y.~Ma, and J.~Wright.
\newblock {Robust Principal Component Analysis?}
\newblock {\em Journal of ACM}, 58(1):1--37, 2009.

\bibitem{cheney2001mathematical}
M.~Cheney.
\newblock {A mathematical tutorial on synthetic aperture radar}.
\newblock {\em SIAM review}, 43(2):301--312, 2001.

\bibitem{Curlander}
John~C. Curlander and Robert~N. McDonough.
\newblock {\em Synthetic Aperture Radar: Systems and Signal Processing}.
\newblock Wiley-Interscience, 1991.

\bibitem{ding2002time}
Y.~Ding and DC~Munson~Jr.
\newblock {Time-frequency methods in SAR imaging of moving targets}.
\newblock In {\em IEEE International Conference on Acoustics, Speech, and
  Signal Processing, 2002. Proceedings.(ICASSP'02)}, volume~3, pages
  2881--2884, 2002.

\bibitem{ding2000analysis}
Y.~Ding, N.~Xue, and DC~Munson~Jr.
\newblock {An analysis of time-frequency methods in SAR imaging of moving
  targets}.
\newblock In {\em Sensor Array and Multichannel Signal Processing Workshop.
  2000. Proceedings of the 2000 IEEE}, pages 221--225, 2000.

\bibitem{ender1993}
J.~Ender.
\newblock {Detectability of slowly moving targets using a multi-channel SAR
  with an along-track antenna array}.
\newblock In {\em Proceedings of SEE/IEE Conference (SARÕ93), Paris}, pages
  19--22, 1993.

\bibitem{fasino}
D.~Fasino.
\newblock {Spectral properties of Toeplitz-plus-Hankel matrices}.
\newblock {\em Calcolo}, 33:87--98, 1996.

\bibitem{tilli}
D.~Fasino and P.~Tilli.
\newblock {Spectral clustering properties of block multilevel Hankel matrices}.
\newblock {\em Linear Algebra and its Applications}, 306:155--163, 2000.

\bibitem{fienup}
J.~R. Fienup.
\newblock {Detecting Moving Targets in SAR Imagery by Focusing}.
\newblock {\em IEEE Transactions on Aerospace and Electronic Systems},
  37(3):794--809, 2001.

\bibitem{friedlander}
B.~Friedlander and B.~Porat.
\newblock {VSAR: a high resolution radar system for detection of moving
  targets}.
\newblock {\em IEE Proc.-Radar, Sonar Navig.}, 144(4):205--218, 1997.

\bibitem{jao}
J.~K. Jao.
\newblock {Theory of Synthetic Aperture Radar Imaging of a Moving Target}.
\newblock {\em IEEE Transactions on Geoscience and Remote Sensing},
  39(9):1984--1992, 2001.

\bibitem{Jakowatz}
Charles V.~Jakowatz Jr., Daniel~E. Wahl, Paul~H. Eichel, Dennis~C. Ghiglia, and
  Paul~A. Thompson.
\newblock {\em Spotlight-mode synthetic aperture radar: A signal processing
  approach}.
\newblock Springer, New York, NY, 1996.

\bibitem{szego}
M.~Kac, W.~L. Murdock, and G.~Szeg\H{o}.
\newblock {On the eigenvalues of certain Hermitian forms}.
\newblock {\em J. Rational Mech. Anal.}, (1):767--800, 1953.

\bibitem{kirscht}
M.~Kirscht.
\newblock {Detection and imaging of arbitrarily moving targets with
  single-channel SAR}.
\newblock {\em IEE Proc.-Radar Sonar Navig.}, 150(1):1984--1992, 2003.

\bibitem{sesana}
E.~Ngondiep, S~Serra-Capizzano, and D.~Sesana.
\newblock {Spectral Features and Asymptotic Properties for $g$-Circulants and
  $g$-Toeplitz Sequences}.
\newblock {\em SIAM J. Matrix Anal. Appl.}, 31(4):1663--1687, 2010.

\bibitem{perry}
R.~P. Perry, R.~C. DiPietro, and R.~L. Fante.
\newblock {SAR Imaging of Moving Targets}.
\newblock {\em IEEE Transactions on Aerospace and Electronic Systems},
  35(1):188--200, 1999.

\bibitem{sparr-time}
T.~Sparr.
\newblock {Time-Frequency Signatures of a Moving Target in SAR Images}.
\newblock {Paper presented at the RTO SET Symposium on Target Identification
  and Recognition Using RF Systems, Oslo, Norway, 11-13 October, 2004.
  Published in RTO-MP-SET-080}.

\bibitem{wang2006}
G.~Wang, X.~Xia, and V.~Chen.
\newblock {Dual-Speed SAR Imaging of Moving Targets}.
\newblock {\em IEEE Transactions on Aerospace and Electronic Systems},
  42(1):368--379, 2006.

\bibitem{wang2004}
G.~Wang, X.~Xia, V.~Chen, and R.~Fiedler.
\newblock {Detection, Location, and Imaging of Fast Moving Targets Using
  Multifrequency Antenna Array SAR}.
\newblock {\em IEEE Transactions on Aerospace and Electronic Systems},
  40(1):345--355, 2004.

\bibitem{zhuMTI}
S.~Zhu, G.~Liao, Y.~Qu, Z.~Zhou, and X.~Liu.
\newblock {Ground Moving Targets Imaging Algorithm for Synthetic Aperture
  Radar}.
\newblock {\em IEEE Transactions on Geoscience and Remote Sensing},
  49(1):462--477, 2011.

\end{thebibliography}
% ---------------------

\end{document}